\documentclass[journal]{IEEEtran}

\ifCLASSINFOpdf
   \usepackage[pdftex]{graphicx}
   \graphicspath{{../figures/}}
   \DeclareGraphicsExtensions{.pdf,.jpeg,.png}
\else
\fi

\ifCLASSOPTIONcompsoc
 \usepackage[caption=false,font=normalsize,labelfont=sf,textfont=sf]{subfig}
\else
 \usepackage[caption=false,font=footnotesize]{subfig}
\fi

\usepackage{url}
\usepackage[colorlinks]{hyperref}
\hypersetup{
	colorlinks=true, 
	linkcolor=blue!70!black, 
	citecolor=red!70!black, 
	filecolor=blue!70!black, 
	urlcolor=blue!70!black, 
}

\RequirePackage{xcolor} 

\definecolor{PairedA}{RGB}{166, 206, 227}
\definecolor{PairedB}{RGB}{ 31, 120, 180}
\definecolor{PairedC}{RGB}{178, 223, 138} 
\definecolor{PairedD}{RGB}{ 41, 128,  35}   
\definecolor{PairedE}{RGB}{251, 154, 153}
\definecolor{PairedF}{RGB}{182,  21,  22}   
\definecolor{PairedG}{RGB}{253, 191, 111}
\definecolor{PairedH}{RGB}{255, 127,   0}
\definecolor{PairedI}{RGB}{202, 178, 214}
\definecolor{PairedJ}{RGB}{106,  61, 154}
\definecolor{PairedK}{RGB}{255, 255, 153}
\definecolor{PairedL}{RGB}{177,  89,  40}

\RequirePackage[cmex10,intlimits]{amsmath}
\RequirePackage{amsfonts}
\RequirePackage{amssymb}

\providecommand{\J}{\ensuremath{\mathrm{j}}}    
\providecommand{\RE}{\ensuremath{\mathrm{Re}}}  

\providecommand{\Quot}[1]{``{#1}"}              
\providecommand{\D}{\,\mathrm{d}}               
\providecommand{\V}[1]{\boldsymbol{#1}}         
\providecommand{\M}[1]{\mathbf{#1}}             
\providecommand{\T}[1]{\mathrm{#1}}             

\providecommand{\herm}{\mathrm{H}} 
\providecommand{\trans}{\mathrm{T}}

\providecommand{\basisFcn}{\V{\psi}}

\providecommand{\XEmat}{\ensuremath{\M{X}_\T{e}}}
\providecommand{\XMmat}{\ensuremath{\M{X}_\T{m}}}

\providecommand{\Prad}{P_\T{rad}}

\providecommand{\We}{\ensuremath{W_\mathrm{e}}}
\providecommand{\Wm}{\ensuremath{W_\mathrm{m}}}

\newcommand{\ie}{\textit{i.e.}{}}
\newcommand{\eg}{\textit{e.g.}{}}

\usepackage{lipsum}
\usepackage{pdfpages}
\usepackage{booktabs}
\usepackage{cancel}
\usepackage[normalem]{ulem}

\providecommand{\figref}[1]{Fig.~\ref{#1}}

\providecommand{\secref}[1]{Section~\ref{#1}}
\providecommand{\appref}[1]{Appendix~\ref{#1}}

\providecommand{\normalvec}{\ensuremath{\V{n}}}  

\providecommand{\des}{\ensuremath{\rho_t}}  
\providecommand{\desD}{\ensuremath{\widetilde{\rho_t}}} 
\providecommand{\desF}{\ensuremath{\bar{\widetilde{\rho_t}}}} 
\providecommand{\Sf}{\ensuremath{S_{f}}} 
\providecommand{\deltarhoMax}{\ensuremath{\Delta \rho_t^\T{max}}} 

\providecommand{\adj}{\ensuremath{\V{\lambda}}} 

\providecommand{\Rmin}{\ensuremath{R_\T{min}}}
\providecommand{\Rs}{\ensuremath{R_\T{s}}}  
\providecommand{\Zmet}{\ensuremath{R_\T{met}}}  
\providecommand{\Zair}{\ensuremath{R_\T{vac}}}  

\providecommand{\densityFilter}{\mathcal{D}}  
\providecommand{\heavisideFilter}{\mathcal{H}}  

\providecommand{\Qe}{\ensuremath{Q_\mathrm{e}}}
\providecommand{\Qm}{\ensuremath{Q_\mathrm{m}}}
\providecommand{\Qem}{\ensuremath{Q_\mathrm{e/m}}}
\newcommand*{\thr}{\raisebox{1.5pt}{\rlap{$\ulcorner$}\kern2pt}} 

\providecommand{\Zvac}{\ensuremath{\M{Z}_0}}  
\providecommand{\Zrho}{\ensuremath{\M{R}_\rho}}  
\providecommand{\Zrhoe}{\ensuremath{\M{\Psi}_\rho}^t}

\providecommand{\Ei}{\ensuremath{\V{E}^\T{i}}}  
\providecommand{\Es}{\ensuremath{\V{E}^\T{s}}}

\begin{document}

\pagestyle{headings}
\twocolumn

\title{Density-Based Topology Optimization in \\ Method of Moments: Q-factor Minimization}

\author{Jonas~Tucek,
        Miloslav~Capek,~\IEEEmembership{Senior~Member,~IEEE,}
        Lukas~Jelinek,
        and~Ole~Sigmund
\thanks{Manuscript received XXX; revised XXX.
This work was supported by the Czech Science Foundation under
project No. 21-19025M, {by the Czech Technical University in Prague
under project SGS22/162/OHK3/3T/13,} and by financial support from the Villum Foundation through the Villum Investigator Project InnoTop.
}
\thanks{J.~Tucek, M.~Capek and L.~Jelinek are with the Czech Technical University in Prague, Prague, Czech Republic (e-mails: {jonas.tucek;~miloslav.capek;~lukas.jelinek}@fel.cvut.cz).}
\thanks{O.~Sigmund is with Technical University of Denmark, Lyngby, Denmark (e-mail: olsi@dtu.dk).}
}


\maketitle

\begin{abstract}
Classical gradient-based density topology optimization is adapted for method-of-moments numerical modeling to design a conductor-based system attaining the minimal antenna  Q-factor evaluated via an energy stored operator. Standard topology optimization features are discussed, \eg{}, interpolation scheme and density and projection filtering. 
The performance of the proposed technique is demonstrated in a few examples in terms of the realized Q-factor values and necessary computational time to obtain a design. The optimized designs are compared to the fundamental bound and well-known empirical structures. The presented framework can provide a completely novel design, as presented in the second example. 
\end{abstract}


\IEEEpeerreviewmaketitle

\section{Introduction}

Wireless communication systems impose more and more requirements on the performance of antenna systems, \eg, compactness, efficiency, or large operational bandwidth. Standard antenna design methodologies~\cite{Balanis_Wiley_2005} are often unsuccessful in delivering an antenna satisfying these demands which are contradictory in nature~\cite{Fujimoto_Morishita_ModernSmallAntennas},~\cite{VolakisChenFujimoto_SmallAntennas}, and designers are forced to utilize advanced procedures based on numerical modeling.

Antenna synthesis attempts to determine the spatial material distribution and feeding network to meet performance requirements and satisfy constraints. This paper focuses solely on optimizing material distribution with predetermined and fixed feeding properties {to reduce the complexity of the optimization problem.}

Antenna engineers often exploit their previous empirical knowledge when fully describing a structure with a set of parameters which are further optimized by simple parametric sweeps or surrogate-based modeling~\cite{koziel2014antenna} in combination with local~\cite{koziel2021robust} or global  optimization methods~\cite{RobinsonSamii_PSOinElmag}. The reliability and  efficiency of these methods are successfully being improved by the community, \eg{}, by achieving immunity to poor initial designs~\cite{koziel2023improved}. Nevertheless, despite reducing the problem's complexity, these techniques struggle to reveal new and non-intuitive designs. Thus, we aim to freely optimize the spatial material presence in a region representing the design domain to eliminate designer bias. The task of optimizing material distribution is called topology optimization~\cite{BendsoeSigmund_TopologyOptimization} which exists in a variety of formulations.

Topology optimization based on genetic algorithms~\cite{RahmatMichielssen_ElectromagneticOptimizationByGenetirAlgorithms} and exact reanalysis~\cite{2021_capeketal_TSGAmemetics_Part1,2021_capeketal_TSGAmemetics_Part2}, attempts to confront the underlying binary nature of the problem~\cite{Nemhauser_etal_IntegerAndCombinatorialOptimization} to find an approximate solution in a large solution space with many local minima. However, these methods are often limited to problems with relatively few design variables due to the necessity of solving a state equation one more time for each design variable for each iteration, resulting in huge amounts of function evaluations~\cite{Sigmund_OnTheUselessOfNongradinetApproachesInTopoOptim}.

In contrast to the aforementioned combinatorial methods, density-based topology optimization~\cite{BendsoeSigmund_TopologyOptimization} applies a continuous relaxation to the binary problem, \ie{}, it replaces binary design variables with continuous ones, and, thus, can exploit adjoint sensitivity analysis~\cite{tortorelli1994design} to compute sensitivities independently on the number of design variables. Consequently, gradient information for objectives is evaluated without significant computational overhead since only one extra adjoint equation needs to be solved per iteration.
Density-based topology optimization is a deterministic local algorithm that can converge within a few hundred function evaluations for even a billion design variables~\cite{Aage_etal_TopologyOptim_Nature2017}. 

The density formulation for topology optimization has been intensely researched during recent decades and has been successfully utilized in various scientific and engineering branches. Apart from the original mechanical problems~\cite{bendsoe1988generating},~\cite{BendsoeSigmund_TopologyOptimization}, it has also been utilized in electromagnetic applications with various numerical techniques, such as the optimization of the distribution of dielectric material, based on finite element modeling~(FEM), \eg{}, to confine light in photonic crystals~\cite{sigmund2008geometric}, and to enhance the transmission properties of nano-photonic waveguides~\cite{jensen2004systematic}, or, based on finite-difference time-domain modeling~(FDTD), \eg{}, to improve the bandwidth of dielectric resonator antennas~\cite{nomura2007structural}. A first approach distributing conductive material was introduced in~\cite{ErentokSigmund2011} and extensively studied in FEM to optimize $s$-parameters of microwave waveguide filters~\cite{aage2017topology} to maximize the radiation efficiency of electrically small antennas~\cite{ErentokSigmund2011}, or, based on FDTD, to minimize return loss of monopole antennas~\cite{HassanWadbroBerggren_TopologyOptimizationOfMetallicAntennas}. 

Density-based topology optimization has also been combined with method of moments (MoM) for the maximization of the total efficiency of an antenna~\cite{wang2017novel}. In that work, the authors did not regularize the small structural features by density filtering~\cite{bruns2001topology}; thus, their designs are mesh-dependent. Furthermore, the fractional bandwidth was not considered as an objective. In~\cite{CismasuGustafsson_AntennaBWoptimizationByGAwithSingleFrequencySimulation}, the bandwidth is optimized through the associated (inversely proportional,~\cite{YaghjianBest_ImpedanceBandwidthAndQOfAntennas}) Q-factor by topology optimization based on genetic algorithms~\cite{RahmatMichielssen_ElectromagneticOptimizationByGenetirAlgorithms} producing highly irregular designs.

This work aims to develop an automated procedure for optimizing parameters of electrically small antennas~\cite{Fujimoto_Morishita_ModernSmallAntennas} within the MoM paradigm~\cite{Harrington_FieldComputationByMoM}. In particular, the Q-factor~\cite{YaghjianBest_ImpedanceBandwidthAndQOfAntennas} is minimized in this paper to demonstrate the efficacy of the procedure.

The projection approach~\cite{wang2011projection} is employed in this paper for gradient-based topology optimization, with the density filtering approach acting as an intermediate step to converge to a similar design with mesh refinement~\cite{bruns2001topology}. The projection scheme is utilized to remedy the introduction of intermediate densities, known as gray scales, which are difficult to interpret in terms of realistic material properties and to ensure convergence to near-binary results. If a conversion from a near-binary to a pure-binary design, consisting only of void and material, is performed, only a slight drop in performance is expected.

The performance of the proposed technique is assessed based on distance to fundamental bounds~\cite{CapekGustafssonSchab_MinimizationOfAntennaQualityFactor}. It is demonstrated on two examples of design domains where human-generated designs exist and possess low values of Q-factor, \ie, a rectangular region and a meanderline~\cite{Capek_etal_2019_OptimalPlanarElectricDipoleAntennas} or a spherical region and a spherical helix~{\cite{Kim_MinimumQESA}},~\cite{Best_LowQelectricallySmallLinearAndEllipticalPolarizedSphericalDipoleAntennas}. The efficiency of these empirical designs can be confirmed by the proposed technique delivering similar topologies or refuted by providing novel better-performing structures.

The paper is structured as follows. The basics of density-based topology optimization are introduced into MoM models in~\secref{sec:TO in MoM}, and the design parametrization is proposed. A suitable interpolation function associating intermediate values of design variables to material parameters is also discussed. Filtering techniques are developed and utilized to penalize unwanted behavior in the designs. A flowchart of the computations of the developed topology optimization framework is summarized in~\secref{sec:flowchart}.
A code for a density-based topology optimization minimizing Q-factor, with annotated comments, is downloadable at~\cite{githubTopOpt}. \secref{sec:examples} proposes a method that minimizes the Q-factor on the two design domains. The method's technicalities are discussed and its time performance is compared to topology optimization by genetic algorithms~\cite{RahmatMichielssen_ElectromagneticOptimizationByGenetirAlgorithms} and to the method based on exact reanalysis~\cite{2021_capeketal_TSGAmemetics_Part1}. The quality of existing empirical designs on a spherical shell is demonstrated by the proposed technique and briefly discussed. The paper is concluded in~\secref{sec:conclusion}.

\section{Topology Optimization in Method-of-Moments Models}\label{sec:TO in MoM}

\begin{figure}[!t]
\centering
\includegraphics[width=3.25in]{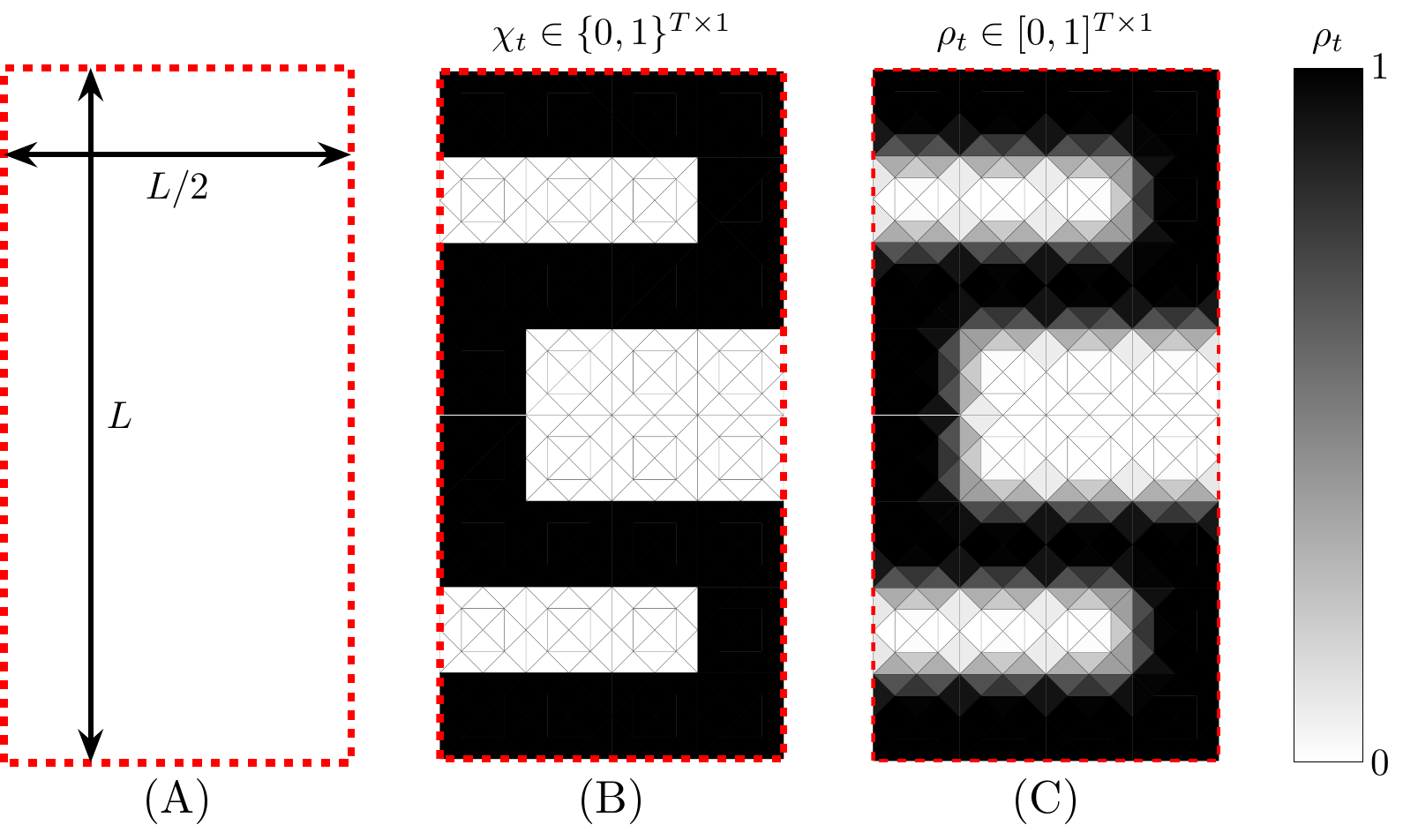}
\caption{(A) Bounding box for the rectangular design domain of~$L\times L/2$ dimensions, (B) design domain discretized into $T$ triangles with a particular binary layout of material represented by material function~$\chi_t$, (C) continuous distribution of material represented by variable~$\des$. Variable~$\des$ represents a density in each element of the grid.
}
\label{fig:Fig1_Discretization}
\end{figure}

Highly conductive surfaces are considered for optimization throughout this paper. The design region is delimited by a bounding box in which the optimized structure is to be located, see the red dashed line in~\figref{fig:Fig1_Discretization}A. The electromagnetic field can be unambiguously represented by an equivalent current density located on the design surface and by the electric field integral equations (EFIE)~\cite{Gibson_MoMinElectromagnetics}, see~\appref{app:EFIE} for further details. The physical problem is made mathematically tractable by discretizing the design domain where unknown state variables will reside. Consequently, current density is expanded into a weighted sum of basis functions and Galerkin's procedure~\cite{Harrington_FieldComputationByMoM} is applied to solve the system equation. Without loss of generality, let us consider triangular discretization and expansion of current density into a set of $N$ Rao-Wilton-Glisson basis functions~\cite{RaoWiltonGlisson_ElectromagneticScatteringBySurfacesOfArbitraryShape}. This leads to a system of linear equations for unknown current expansion coefficients~$\M{I}\in\mathbb{C}^{N\times 1}$ reading
\begin{equation}
    \M{ZI}=\left(\Zvac + \Zrho \right)\M{I} = \M{V},
    \label{eq:ZI=V}
\end{equation}
where $\Zvac\in\mathbb{C}^{N\times N}$ and $\Zrho \in \mathbb{R}^{N\times N}$ are the vacuum and material parts of the impedance matrix, respectively. Vector~$\M{V}\in \mathbb{C}^{N\times 1}$ represents excitation coefficients, \eg, a plane-wave or a delta gap feeding. Throughout this paper, the excitation vector $\M{V}$ is given \textit{a priori} and is, therefore, fixed during the optimization. The optimization of volumetric dielectric bodies can be equally defined with minimum modifications employing a volumetric electric field integral equation, however, such an approach is out of the scope of this paper.

For a given excitation~$\M{V}$, the system~\eqref{eq:ZI=V} defines a physically admissible response (current density)
\begin{equation}
    \M{I} = (\Zvac + \Zrho)^{-1}\M{V},
\end{equation}
and, therefore, formula~\eqref{eq:ZI=V} is a fundamental optimization constraint.  

\subsection{Design Parametrization}

In this paper, the goal of topology optimization is to distribute a highly conducting material, a perfect electric conductor (PEC), over a fixed design domain to meet the requirements of antenna performance. To facilitate optimization, the design domain is the same as the modeling domain, and each triangular element, referred to as a \Quot{pixel}, is parametrized with surface resistivity~$\Rs$. Hence, the topology optimization task is reduced to a classical material distribution task where we aim to determine a material characteristic function 
\begin{equation}
\chi_t (\V{r}) = 
    \begin{cases}
    0\quad\iff\quad \Rs\approx \infty\;\Omega~\T{(vacuum)}, \\
    1\quad\iff\quad \Rs\approx 0\;\Omega~\T{(PEC)},
    \end{cases}
\end{equation}
where $\chi_t$ controls the layout of two materials: air and a PEC, see \figref{fig:Fig1_Discretization}B. 

Density-based topology optimization~\cite{BendsoeSigmund_TopologyOptimization} commonly replaces the discrete material function~$\chi_t$ with the continuous design variable as
\begin{equation}
    0 \leq \des \leq 1,
\end{equation}
so that the design variable is allowed to take intermediate values between the air and the PEC representing a \Quot{density} of material, see \figref{fig:Fig1_Discretization}C. Subscript $_t$ represents design variable $\des$ as a piece-wise constant over the triangular domain given by the~$t$-th triangle. 

Consequently, the interpolation function~\cite{bendsoe1999material} is used to associate surface resistivity~$\Rs$ with intermediate values of~$\des$ and the physics is projected to material matrix~$\Zrho$ as
\begin{equation}
    \Zrho = \displaystyle\sum_{t=1}^T \Rs(\des)\Zrhoe,
\end{equation}
where each pixel attains real constant surface resistivity~$\Rs$ parametrized by design variable $\des$ through an interpolation function, and $\Zrhoe$ is a material element matrix, see~\appref{app:EFIE} for its form.

Generating a suitable type of interpolation function is, generally, the most crucial step within topology optimization algorithms~\cite{bendsoe1999material}. The function should vary monotonically, and a slight change in~$\des$ should produce a non-negligible response in function value~$\Rs$. Surface resistivity acts as a damping parameter for electromagnetic field values, \ie{}, its magnitude has a strong effect, and thus, the range of resistivities in the design plays a crucial role, both for the physics evaluation and optimization process. Consequently, we cannot adapt the simple linear interpolation function, which has been successfully employed in dielectric designs~\cite{kiziltas2003topology},~\cite{jensen2004systematic},~\cite{nomura2007structural}, mainly because of the huge contrast between the surface resistivity of vacuum and PEC. Therefore, a more suitable interpolation scheme for our design approach is developed in Section~\ref{sec:InterpScheme}. 

\subsection{Interpolation Scheme}\label{sec:InterpScheme}
Considering the aforementioned requirements on the interpolation function, it is appropriate to express surface resistivity according to~\cite{ErentokSigmund2011} and~\cite{kiziltas2003topology} as
\begin{equation}
    \Rs(\des) = \Zair \left(\dfrac{\Zmet}{\Zair} \right)^{f(\des)}, 
    \label{eq:Interpolation}
\end{equation}
where $f(\des)$ is a penalization function based on a RAMP-like\footnote{Rational Approximation of Material Properties.} model~\cite{BendsoeSigmund_TopologyOptimization}
\begin{equation}
    f(\des) = \frac{\des}{1+p(1-\des)},
\end{equation}
where $p$ is a penalty parameter used to regulate the behavior of the interpolation and penalize intermediate values. In this work, $p=1$, and we utilize other penalization techniques shown later in~\secref{sec:Filtering}. Surface resistivity~\eqref{eq:Interpolation} attains range $[\Zair,\Zmet]$, whose boundaries represent the resistivity of the vacuum and PEC.

The upper bound~$\Zmet=\Rs(\des=1)$ and the lower bound~$\Zair=\Rs(\des=0)$ on attainable resistivities should be set to simulate the actual physical response of the material distribution correctly but, at the same time, have the smallest possible contrast. They are set based on the numerical investigations of the authors, see~\appref{app:Choosing resistivities}, as 
\begin{equation}
    \Zair = 10^5\;\Omega, \quad \Zmet=1\;\Omega,
    \label{eq:resistivities bounds}
\end{equation}
which, for this paper (Q-factor minimization), provides physically acceptable results and enough freedom for the optimization step since the Q-factor is not influenced by loss power, see~\appref{app:Choosing resistivities} for more details. Therefore, resistivity bounds~\eqref{eq:resistivities bounds} represent the approximative physical behavior of the vacuum and PEC, but a subsequent analysis of the final design using the true resistivities, \ie{},~$\Zair\approx\infty\;\Omega$ and~$\Zmet \approx 0\;\Omega$, results in a negligible drop in Q-factor value. Nevertheless, current density is influenced, yet it is important to emphasize that these values are not universal and cannot be used for arbitrary physical metrics. Numerical investigations should always be conducted to determine the properties for a particular optimization goal.  

\begin{figure}[!t]
\centering
\includegraphics[width=3.25in]{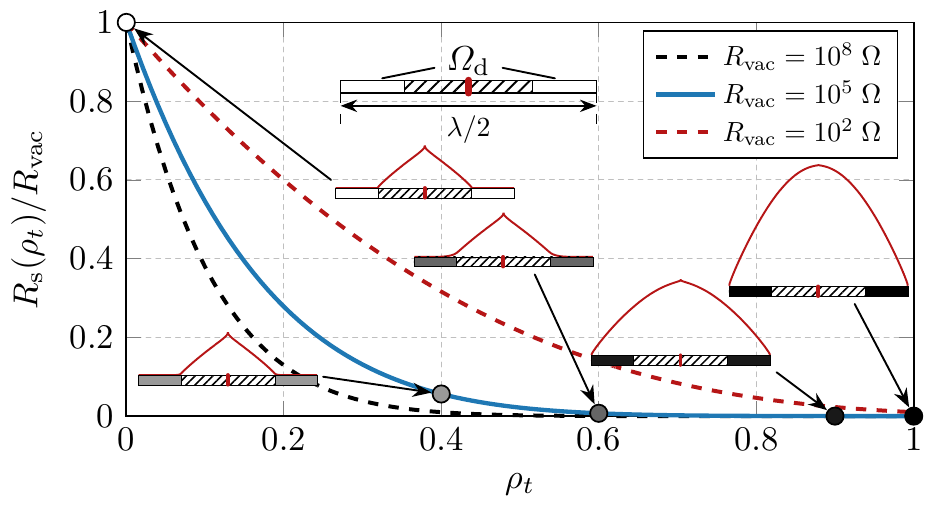}
\caption{Interpolation scheme of surface resistivity $\Rs(\des)$ as a function of density~$\des$ for three different values of vacuum resistivity~$\Zair$. The influence of the interpolation on a current density is explored on a half-wave-long strip dipole antenna with the ends being designable domain~$\varOmega_\T{d}$, in which $\des$ is varied from 0 to 1. The Red solid lines show a current distribution along the dipole excited by a delta gap placed in the middle, \eg{}, a standard $\lambda/2$-long dipole current is observed for $\des=1$.}
\label{fig:Fig2_Interpolation}
\end{figure}

The interpolation function~\eqref{eq:Interpolation} is shown in~\figref{fig:Fig2_Interpolation} for three values of $\Zair$. The exponential dependency on~$\des$ results in a fast decay from $\des=1$, which is demonstrated on a current distribution on a center-fed dipole antenna in which the end parts are covered by variable resistivity ($\des$ is varied from 0 to 1) while the central part is covered with resistivity~$\Zmet=\Rs(\des=1)$. We observe standard cosinusoidal and triangular current distribution for the half-wave and a short dipole, respectively. 

\begin{figure}[h!]
\centering
\includegraphics[width=3.25in]{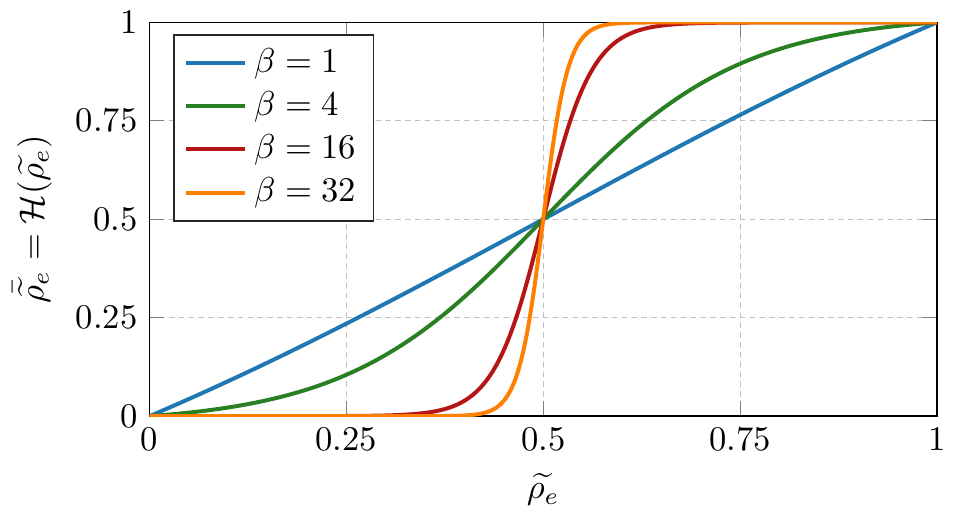}
\caption{{Projection filter~$\heavisideFilter$ with four different values of sharpness~$\beta$ and projection level~$\eta=0.5$.}}
\label{fig:Fig3_Heaviside}
\end{figure}

\begin{figure*}[h!]
\centering
\includegraphics[width=1\textwidth]{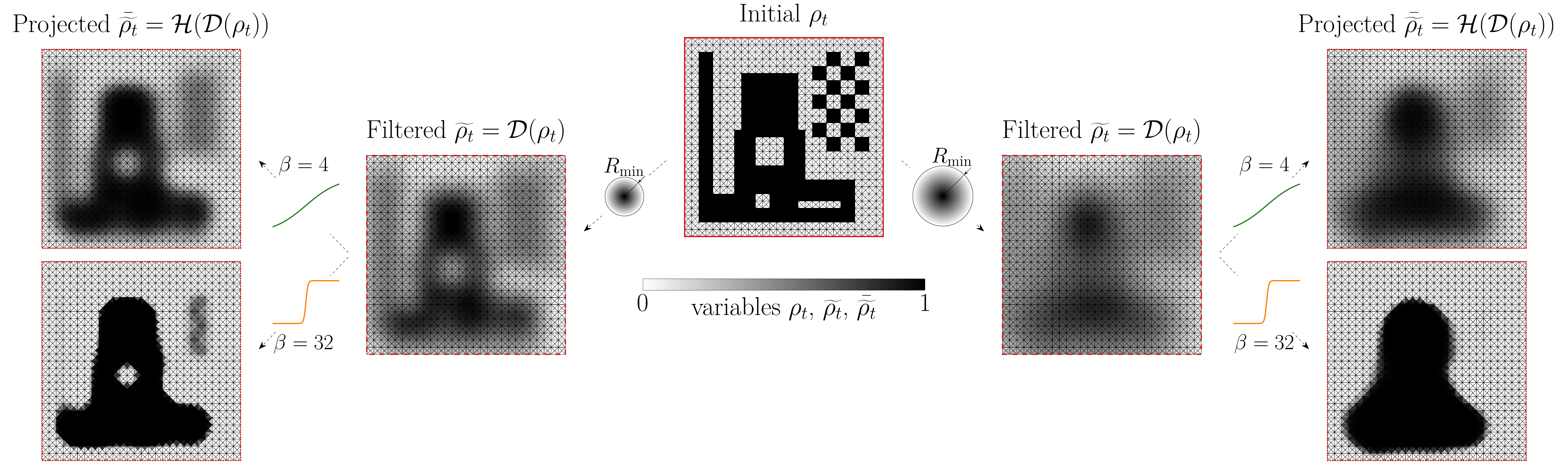}
\caption{ {Demonstration of the filtering process.} The {initial} test design {(solid box)} represented by $\des$ contains specific patterns, \eg, sheets, and holes with various widths or checkerboard patterns, comparable in size with density filter radius~$\Rmin$. Density filtering of these patterns results in distribution~$\desD$ containing gray areas {(dashed boxes)}, possibly removing undesirable details by smoothing the design distribution, thus blurring the solution space which might result in an escape from the local minimum. The subsequent application of the projection filter removes unwanted gray areas resulting in a design of near-binary nature {(dotted boxes)}.}
\label{fig:Fig4_Filters}
\end{figure*}

\subsection{Minimization of Q-Factor}
For single resonance systems, Q-factor provides a measure of the frequency selectivity of an antenna~\cite{Fante_QFactorOfGeneralIdeaAntennas}, \ie, an estimate of a fractional bandwidth through its inverse proportionality $B\sim Q^{-1}$~\cite{YaghjianBest_ImpedanceBandwidthAndQOfAntennas}, see~\appref{app:Quality factor} for its definition. For electrically small antennas the lowest Q-factor is observed at self-resonance~\cite{CapekGustafssonSchab_MinimizationOfAntennaQualityFactor},~\cite{Chu_PhysicalLimitationsOfOmniDirectAntennas}. Nevertheless, previous design experience is required to construct such a structure, \eg{}, a spherical helix antenna~\cite{Best_LowQelectricallySmallLinearAndEllipticalPolarizedSphericalDipoleAntennas}, or a meanderline antenna~\cite{Capek_etal_2019_OptimalPlanarElectricDipoleAntennas}.

The density-based topology optimization is carried out in this paper to find an optimized layout of conductive material to minimize the antenna Q-factor within MoM models. The optimization formulation is given by
\begin{equation}
    \begin{aligned}
        \underset{\des}{\T{minimize}} &\quad Q(\M{I})=\T{max}\{\Qe(\M{I}),\Qm(\M{I})\} \\
         \T{subject~to} 
         &\quad \left[\Zvac + \Zrho(\des)\right]\M{I} = \M{V},\\
         &\quad \dfrac{1}{A_0} \sum_{t=1}^T\des A_t\leq \Sf,\\
         &\quad 0\leq \des \leq 1,\quad t=1,\dotsc,T,
         \label{eq:minmax}
    \end{aligned}
\end{equation}
where subscripts $_{\T{e}/\T{m}}$ denote the electric and magnetic Q-factor parts, respectively, $\des$ is the density design variable associated with each triangular pixel, $A_t$ is the area of the $t$-th triangle, $A_0$ is the design domain area, and $\Sf$ is the maximum fraction of design domain area which can be occupied with material. Task~\eqref{eq:minmax} is defined as a nested formulation, \ie, the MoM equation is included for the sake of completeness but its solution is directly incorporated into the evaluation of the objective and constraints~\cite{BendsoeSigmund_TopologyOptimization}. Although electromagnetic designs do not benefit from spanning the whole design domain with material, the second constraint, which restricts the amount of the available material, is applied to speed up the convergence and to limit the existence of non-connected regions of the conductor.

\subsection{Filtering Techniques}\label{sec:Filtering}

Solving topology optimization task~\eqref{eq:minmax} as it stands often results in mesh-dependent designs. Whenever a mesh is refined, the topology optimization may add details comparable to the mesh element size and often converge to a significantly different topology. Moreover, checkerboard-like patterns might occur in this resulting design~\cite{diaz1995checkerboard}. The continuous relaxation also inevitably leads to the appearance of residual gray pixels, \ie{}, intermediate densities, in the actual design. This results in final designs being practically unrealizable. 

A projection operator~\cite{wang2011projection} is used to project continuous design variables into a binary space as
\begin{equation}
    {\desF} = \heavisideFilter({{\desD}}),
\end{equation}
where {$\desF$} is the projected {filtered} density and $\heavisideFilter$ is the Heaviside function. The Heaviside function does not have a finite derivative and is thus replaced by a continuous function~\cite{wang2011projection} in the form
\begin{equation}
    \heavisideFilter ({\desD}) \approx \frac{\tanh(\beta \eta) + \tanh(\beta({\desD}- \eta))}{\tanh(\beta \eta) + \tanh(\beta(1- \eta))},
    \label{eq:Heaviside}
\end{equation}
where parameters $\beta$ and $\eta$ are threshold sharpness and threshold level, respectively, see~\figref{fig:Fig3_Heaviside}. 

The {original design variables~$\des$ are} smoothed to regularize the solution space of the optimization problem to provide an existence of solutions. Several regularization schemes are reported in the literature, \eg{}, perimeter control~\cite{haber1996new}, sensitivity filtering~\cite{sigmund1997design}, and density filtering~\cite{bruns2001topology},~\cite{bourdin2001filters}.
In this paper, density filtering operation~$\densityFilter(\cdot)$ is used as an intermediate step written in the form of a convolution~\cite{bourdin2001filters} as
\begin{equation}
    \desD = \densityFilter(\des) = \sum_{j\in B_{t}}W_{tj}\rho_j,
    \label{eq:density}
\end{equation}
where $\desD$ is the filtered variable, $\rho_j$ is the original design variable and $W_{tj}$ is a convolution kernel defined as 
\begin{equation}
    W_{tj} = \frac{(\Rmin - \|\V{r}_t-\V{r}_j\|)}{\sum\limits_{{k} \in B_{{t}}}(\Rmin - \|\V{r}_t-\V{r}_{k}\|)},
\end{equation}
with~$\V{r}_m$ being the center of the~$m$-th triangle.
Support domain~$B_{t}$ for the filter is a mesh-independent neighborhood of $t$-th triangle specified by a sphere with radius~$\Rmin$ as
\begin{equation}
    B_{t} \in \left\{j\;:\;\|\V{r}_t-\V{r}_j\| \leq \Rmin\right\}.
\end{equation}

The density filtering is introduced to impose a length-scale, \ie{}, a characteristic length used to control the size of features in the optimized design. It can be interpreted as the minimum feature size the optimizer can generate. The density filter is utilized to regularize the solution space and accelerate the optimization process but introduces intermediate densities {, motivating the subsequent} projection to{wards} binary design~{by~\eqref{eq:Heaviside}}.

The projection operator~\eqref{eq:Heaviside} is utilized in a continuation scheme to converge to a near-binary
design {in order} to moderate adverse effects. Sharpness~$\beta$ is gradually increased during the optimization to minimize the effect of gray elements and to converge to a near-binary result. The pure-binary material distribution is obtained by the hard thresholding so that~$\desF\geq0.5 \to 1$ and~$\desF<0.5 \to 0$, \ie{}, applying~\eqref{eq:Heaviside} with sharpness~$\beta\to\infty$.
 
The influence of filters is demonstrated in~\figref{fig:Fig4_Filters}, where a test design with patterns of varying complexity is filtered by a density filter with two different radii~$\Rmin$ comparable with the included details. Despite the introduction of gray areas by the density filter,  unwanted details, \eg{}, the checkerboard pattern, are removed. The density filtering technique provides the existence of the solution and convergence with the mesh refinement. The subsequent utilization of the projection filter~\eqref{eq:Heaviside} limits the gray areas and ensures convergence to near-binary results, but the length-scale imposed by the density filter is lost. To ensure length-scale control as well, a robust formulation of topology optimization~\cite{wang2011projection} may be exploited, but such a formulation is out of the scope of this work.    

\section{Q-factor Minimization via Topology Optimization}

The optimization task~\eqref{eq:minmax} is regularized by the filtering techniques, \ie{}, density, and projection filters, and requires continuous derivatives to be solved by the gradient-based optimizer. {In order to solve the originally non-differentiable}  \Quot{min-max} problem~\eqref{eq:minmax}{, the problem is converted to the so-called bound formulation~\cite{BendsoeSigmund_TopologyOptimization} }
\begin{equation}
    \begin{aligned}
        \underset{\des}{\T{minimize}} &\quad z \\
         \T{subject~to} 
         &\quad z \geq \Qe(\M{I}),\\
         &\quad z \geq \Qm(\M{I}),\\
         &\quad \left[\Zvac + \Zrho(\desF)\right]\M{I} = \M{V},\\
         &\quad \dfrac{1}{A_0}\sum_{t=1}^T\desF A_t \leq \Sf,\\
         &\quad \desF = \heavisideFilter\left(\densityFilter\left(\des\right)\right),\\
         &\quad 0\leq \des \leq 1,\quad t=1,\dotsc,T,\\
         \label{eq:MMA task}
    \end{aligned}
\end{equation}
where an extra decision variable~$z$ is introduced. We caution{, however,} that the optimization task~\eqref{eq:MMA task} is continuous but non-convex, thus global optimality is not guaranteed. Therefore, the quality of the solution can be solely judged based on the comparison to fundamental bounds~\cite{CapekGustafssonSchab_MinimizationOfAntennaQualityFactor} {or comparison to existing designs}. 

The topology optimization task~\eqref{eq:MMA task} is solved for local optimality using the Method of Moving Asymptotes (MMA)~\cite{svanberg1987method}, which is a gradient-based optimizer requiring knowledge of the sensitivities of objectives and constraints with respect to the design variable{. We here use that all components of objectives and constraints are continuously differentiable}. The sensitivities of $\Qe$ and $\Qm$, with respect to physical design variable~$\desF$, are evaluated by an adjoint sensitivity analysis~\cite{tortorelli1994design} as
\begin{equation}
        \frac{\D \Qem}{\D \desF} = 2 \RE \left\{\adj^\trans_\T{e/m}\frac{\partial \Zrho}{\partial \desF}\M{I} \right\},
        \label{eq:Sensitivities}
\end{equation}
where adjoint vectors $\adj_\T{e/m}$ are determined from the corresponding adjoint equations
\begin{equation}
    \M{Z}^\trans\adj_\T{e/m} = -\left(\frac{\partial \Qem}{\partial \M{I}}\right)^\trans.
    \label{eq:Adjoints}
\end{equation}
 The derivative of $\Zrho$ with respect to $\desF$ is 
\begin{equation}
    \frac{\partial \Zrho}{\partial \desF} =  \frac{\partial R_\T{s}}{\partial \desF}\Zrhoe,
\end{equation}
where $\partial R_\T{s}/\partial\desF$ is the derivative of the interpolation function~\eqref{eq:Interpolation}. Thus, to compute the sensitivity of the Q-factor, only the derivatives $\partial \Qe/\partial\M{I}$ and $\partial \Qm/\partial\M{I}$, see~\appref{app:Quality factor}, need to be evaluated. The system matrix is reused to solve adjoint equations~\eqref{eq:Adjoints} with the corresponding right-hand side.

Sensitivities~\eqref{eq:Sensitivities} are evaluated with respect to the filtered {and projected}, \ie, physical, density field~$\desF$, hence, the chain rule is utilized to determine sensitivities with respect to a change in the design field~$\des$ as   
\begin{equation}
    \frac{\D \Qem}{\D \des} = \sum_{i\in B_{t}} \frac{\D \Qem}{\D \bar{\widetilde{\rho_i}}} \frac{\partial \bar{\widetilde{\rho_i}}}{\partial \widetilde{\rho_i}} \frac{\partial \widetilde{\rho_i}}{\partial \des},
\end{equation}
where the derivative of the projected density~$\desF$, with respect to a change in the filtered density~$\desD$, is determined as
\begin{equation}
    \frac{\partial \bar{\widetilde{\rho_i}}}{\partial \widetilde{\rho_i}} =\frac{ \beta (1-\tanh(\beta(\widetilde{\rho_i} - \eta))^2}{\tanh(\beta \eta) + \tanh(\beta(1- \eta))},
\end{equation}
and where the sensitivity of the filtered field~$\desD$, with respect to design variable~$\des$, is found as
\begin{equation}
    \frac{\partial \widetilde{\rho_i}}{\partial \des} = \frac{(\Rmin - \|\V{r}_t-\V{r}_i\|)}{\sum\limits_{{k} \in B_{t}}(\Rmin - \|\V{r}_t-\V{r}_{{k}}\|)}.
\end{equation}

\begin{figure}[!t]
\centering
\includegraphics[width=3.25in,trim={1.15cm 0 0 0},clip]{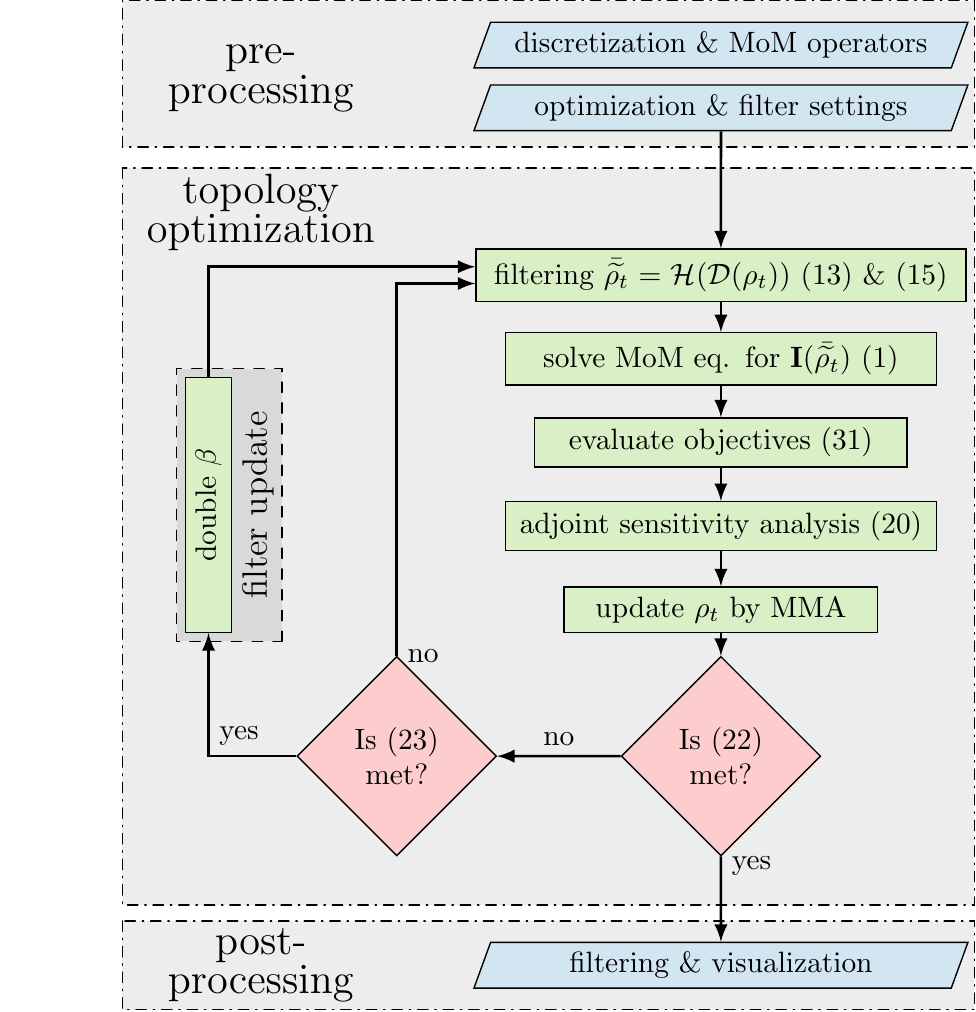}
\caption{Flowchart of the proposed topology optimization procedure.}
\label{fig:Fig5_flowchart}
\end{figure}

\section{Flowchart of the Topology Optimization procedure}\label{sec:flowchart}
To implement the optimization procedure, users can refer to the code implementation of density topology optimization which is available as supplementary material~\cite{githubTopOpt} and includes annotated comments for each part of the process. The flowchart in~\figref{fig:Fig5_flowchart} provides an overview of the optimization stages and is further detailed in this section. During the preprocessing, MoM operators are calculated on a fixed grid, and the user input parameters, such as area fraction ($\Sf$), density filter radius ($\Rmin$), initial design variable ($\des$), maximum iterations ($I$), and maximum design variable change ($\deltarhoMax$, defaulted to 0.01) are set. The projection filter~\eqref{eq:Heaviside} starts at level~$\eta=0.5$ and sharpness~$\beta=1$, which is later updated in the optimization loop until it reaches a maximum of $\beta_\T{max}=32$.

The iterative optimization process follows: filtering, solving the MoM equation, evaluating objectives, computing sensitivities, and updating the design variable~($\des$) by MMA {with move limits~\cite{svanberg1987method} set to 0.25}. If the maximum number of iterations is not reached and the maximum change in~$\des$ in two consecutive iterations~$i$ is greater than $\deltarhoMax$, \ie{}, if conditions 
\begin{equation}
    \underset{{t}}{\T{max}}\big\{|\des^i - \des^{i-1}|\big\} \geq \deltarhoMax \quad\T{and} \quad i < I
\end{equation}
are met, the optimization continues, and a decision about the update of sharpness~$\beta$ of the projection filter is made. The sharpness is doubled if the change in~$\des$ in two consecutive iterations~$i$ is less than~$\deltarhoMax$ or every 100 iterations until the sharpness reaches maximum value, \ie{}, if logic relation
\begin{equation}
    \begin{aligned}
        \big(\underset{{t}}{\T{max}}\big\{|\des^i - \des^{i-1}|\big\}\leq \deltarhoMax \quad&\T{or} \quad \T{mod}(i,100)=1\big)\\
        &\T{and} \quad \beta \leq \beta_\T{max}
    \end{aligned}
\end{equation}
are fulfilled, and where~$\T{mod}()$ is the modulo operation. After meeting the stopping criterion, the design and physical quantities are post-processed and visualized. The continuation scheme for the projection filter produces a near-binary optimized design which is finally converted to a pure-binary design through hard thresholding. The Q-factor of the thresholded design is denoted as $\thr Q$ and is evaluated with physically accurate values of the resistivities for the vacuum~($\Zair\approx\infty\;\Omega$) and PEC~($\Zmet\approx0\;\Omega$).

Optimization parameters significantly impact the quality of the solutions, and multiple optimization runs with different initial guesses or parameters may be necessary to reach the best possible optimized design as the problem is non-convex and can become trapped in local minima. Care must also be taken in selecting $\Rmin$ and the continuation scheme's speed, as a large filter radius can prevent good minima from being reached, and a rapid increase in $\beta$ can destabilize convergence. Testing different parameters on a coarse mesh before running the final optimization with a fine grid is highly recommended.

\begin{figure}[!t]
\centering
\includegraphics[width=3.25in]{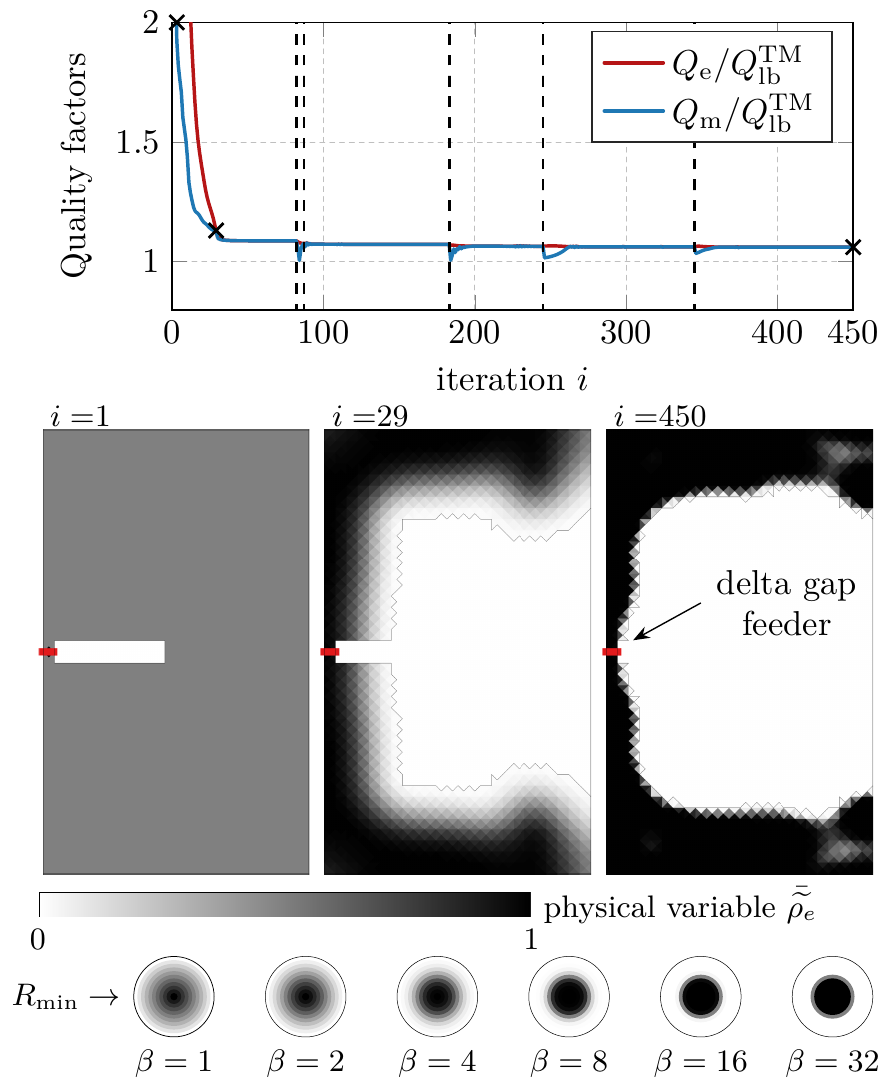}
\caption{(Top) Convergence history of Q-factors $Q=\T{max}\{\Qe, \Qm\}$ normalized to the fundamental bound~$Q_\T{lb}^{\T{TM}}$~\cite{GustafssonTayliEhrenborgEtAl_AntennaCurrentOptimizationUsingMatlabAndCVX}, where dashed vertical lines represent iterations in which the sharpness~$\beta$ of projection filter~$\heavisideFilter$ was doubled. (Bottom) Snapshots of density field distribution from a few selected iterations are included below in the convergence graph. The optimizer reaches a self-resonant design in iteration~$i=29$. The structure is only slightly refined during the rest of the iterations, mainly by removing gray elements by the continuation scheme of the projection filter. The black line bounds the region where~$\desF\geq0.01$. The density filter size characterized by radius~$\Rmin$ is also depicted in scale with the consequent utilization of the projection filter with sharpness~$\beta$, effectively reducing its size and losing the imposed length-scale.}
\label{fig:Convergence_Plate_40x24}
\end{figure}

\section{Numerical Investigations}\label{sec:examples}
The proposed topology optimization technique is implemented in MATLAB~\cite{matlab} according to~\secref{sec:flowchart}, and its properties and performance are discussed. The goal is to minimize the Q-factor~\eqref{eq:MMA task} by distributing a conductive material over a fixed triangular discretized domain. The required MoM matrices are computed in AToM~\cite{atom}. The examples are evaluated on a computer with Intel Xeon Gold 6244 CPU (16 physical cores, 3.6 GHz) with 384\,GB RAM. The excitation of the antennas is performed via delta gap feeders with their immediate neighborhood kept fixed\footnote{The excitation must be connected at all times, and the action of the density filter may cause unwanted disconnections or shortages of the feed.}.

\subsection{Rectangular Design Region}\label{sec:rectangle}

The first example considers a rectangular design region of aspect ratio~$3:5$ and electrical size $ka=0.8$, where~$k$ is the wavenumber and~$a$ is the radius of the smallest circumscribing sphere. A delta gap feeder provides excitation close to the left boundary, see the design at the bottom of~\figref{fig:Convergence_Plate_40x24}. The objectives, \ie{}, Q-factors~$\Qem$, are normalized to the bound $Q_\T{lb}^\T{TM}$~\cite{GustafssonTayliEhrenborgEtAl_AntennaCurrentOptimizationUsingMatlabAndCVX} which is 
valid for antennas radiating TM modes only. The design domain is discretized into 3\,840 triangles (5\,696~basis functions). The Q-factor minimization~\eqref{eq:MMA task} is performed by the density-based topology optimization with a fixed radius of the density filter set to~$\Rmin = 0.15a$, gradually increased sharpness of the projection filter according to~\secref{sec:flowchart}, and the maximal number of iterations $I=600$. The optimization is subjected to area fraction constraint~$\Sf=0.35$, and all design variables~$\des$ are initially set to the prescribed area fraction. 

The optimization ran for 135 minutes, and its progress is reported at the top of~\figref{fig:Convergence_Plate_40x24}. It can be observed that the proposed technique quickly locates a self-resonant structure that separates the electric charge well, \ie{}, an expected solution of the minimal Q-factor~\cite{GustafssonSohlKristensson_IllustrationsOfNewPhysicalBoundOnLinearlyPolAntennas}. The subsequent convergence rate is mostly flat, with occasional jumps resulting from the projection filter update. The design is only slightly refined by further removing the gray elements. The optimized design, see $i=450$ at the bottom of~\figref{fig:Convergence_Plate_40x24}, exhibits a Q-factor very close to the fundamental bound~$Q/Q_\T{lb}^\T{TM}=1.06$. Below the snapshots in~\figref{fig:Convergence_Plate_40x24}, you will find information about the impact of the projection filter on the density filter size while updating sharpness~$\beta$. While the projection scheme yields a near-binary design, it also disrupts the length-scale imposed by the underlying density filter. If no special precautions are taken, the length-scale near the outer boundary may be violated~\cite{clausen2017filter}. The exact length-scale in this paper is {not crucial as long as the final design avoids unwanted features, \eg{}, checkerboard-like patterns~\cite{diaz1995checkerboard}.}{The primary purpose of filtering is} to regularize the solution space and accelerate the {optimization process}.

\begin{figure}[!t]
\centering
\includegraphics[width=3.25in]{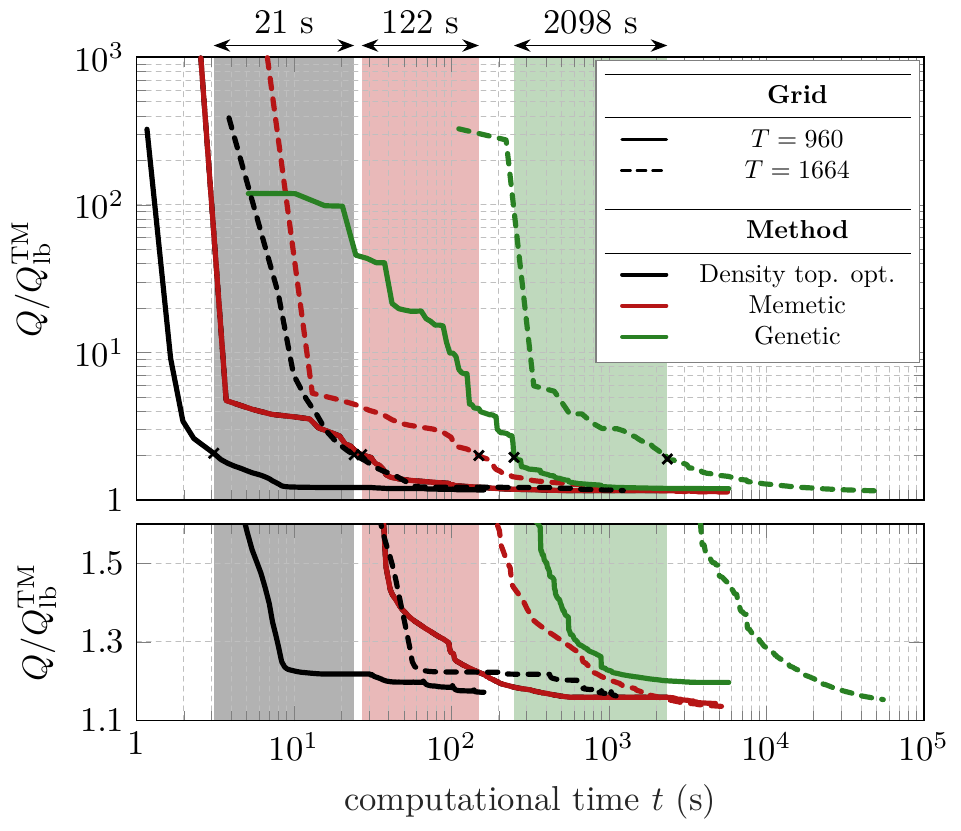}
\caption{Performance of topology optimization algorithms for Q-factor minimization for two grids. Settings are the same as in~\figref{fig:Convergence_Plate_40x24}. The objective is normalized to the fundamental bound~$Q_\T{lb}^{\T{TM}}$. Three implementations are shown: density topology optimization (black curves), a genetic algorithm (green curves) and a memetic combination of global and local step (red curves). Black crosses mark the computational time needed to reach $Q = 2Q_\T{lb}^{\T{TM}}$ and show the rough estimate on the scalability of the method with the number of optimization variables.  
}
\label{fig:Computational complexity}
\end{figure}


{Despite the proposed method being an inherently local procedure, it is resilient towards initial design choices due to the nature of the Q-factor~\cite{GustafssonCismasuJonsson_PhysicalBoundsAndOptimalCurrentsOnAntennas} and the usage of spatial filtering, see~\appref{app:initial seed} for more details. While optimizing the Q-factor, the method delivers locally extremal designs that are in the performance close to the fundamental bound, thus approaching a global extremum.}

The computational time necessary to obtain a design from the proposed optimization is further judged for two discretization grids, see~\figref{fig:Computational complexity}. It is observed that the density topology optimization provides convergence to local minimum and has a steep convergence rate which is typical for gradient-based algorithms. It provides gray designs with low values of the Q-factor in a short computational time. Topology optimization using density-based approaches exhibits good scalability with increasing design variables, as noted in \cite{Aage_etal_TopologyOptim_Nature2017}. For a rough estimate, please refer to the gray region in~\figref{fig:Computational complexity}. The relatively low computational cost allows rerunning the optimization several times for various preferences of user-defined parameters, \eg{}, interpolation function, area fraction, etc., to determine which settings produce the best result.  

 The computation time is also compared with other topology optimization techniques, namely, binary topology optimization based on a genetic algorithm and topology optimization based on memetic scheme~\cite{2021_capeketal_TSGAmemetics_Part1} combining exact reanalysis and a genetic algorithm, see~\figref{fig:Computational complexity}. Both methods utilize 144 agents. For the sake of simplicity, the computation time is not normalized to the number of physical cores used for the evaluation. The convergence rate of the memetic scheme is superior to the sole use of the genetic algorithm, but it is much slower than the density-based approach. An estimate on scalability is also included, see the red and green regions in~\figref{fig:Computational complexity}, and is inferior to the density topology optimization.

\begin{figure}[!t]
\centering
\includegraphics[width=3.25in]{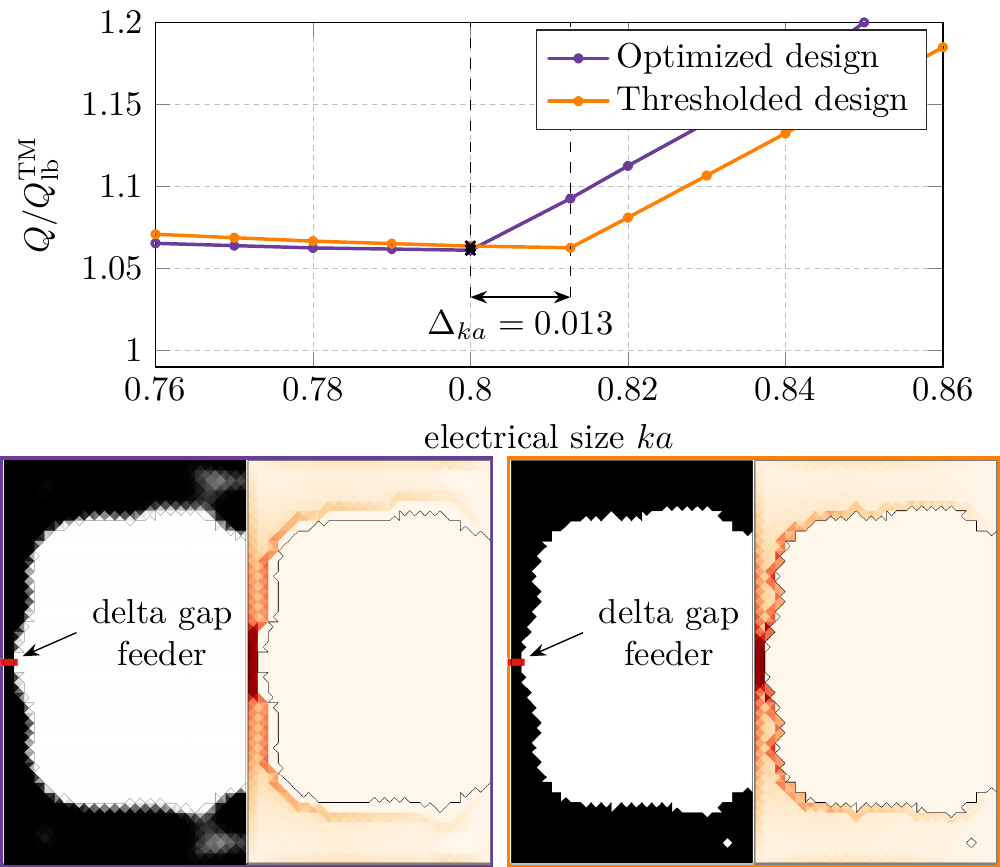}
\caption{(Top) Frequency sweep of the normalized Q-factor of the design optimized at single electrical size $ka=0.8$. The purple curve corresponds to density-based topology optimization with gradually increased sharpness of the projection filter. The orange curve corresponds to the thresholded design. The thresholding may result in a slight frequency shift and performance drop. (Bottom left) Optimized and (bottom right) thresholded designs are also included with the dipole-like current density excited by a delta gap feeder depicted by the blue line. The colormap shows an absolute value of the current density. The black line bounds the region where~$\desF\geq0.01$.}
\label{fig:Plate_40x24_kaSweep}
\end{figure}

\begin{figure*}[!t]
\centering
\includegraphics[width=1\textwidth]{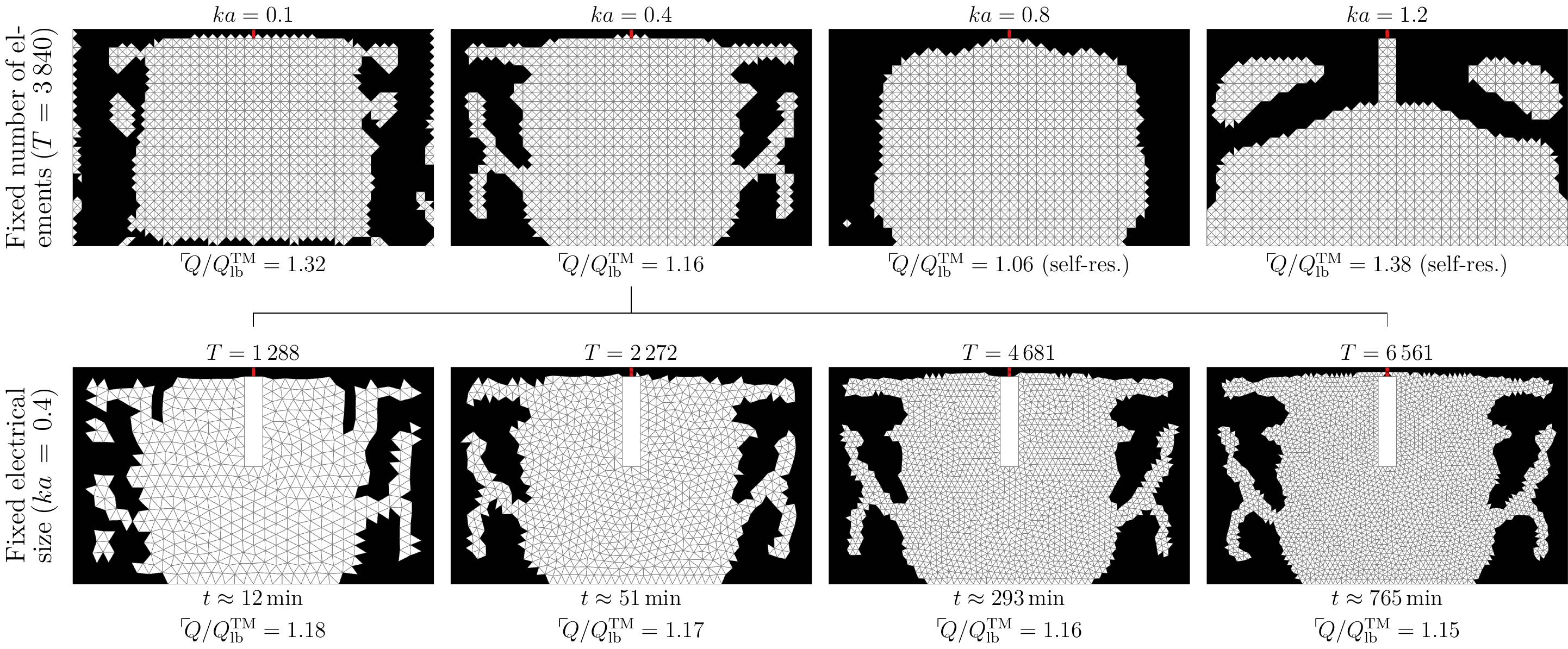}
\caption{The proposed optimization technique is utilized to minimize the Q-factor at the various electric sizes and fixed discretization~($T=3840$). The first row shows the thresholded design and corresponding Q-factor. The grid cannot accommodate a self-resonant structure for the small electric sizes, and much finer granularity would be required. In the second row, the optimization is repeated for a particular electric size~$ka=0.4$ while the discretization grid is changed to an unstructured one and further refined. The density filter with a fixed radius provides the solution while the mesh is refined.}
\label{fig:Other Designs}
\end{figure*}

Despite the use of a continuation scheme, the optimized design obtained by the density topology optimization still contains residual gray elements, see the design in the bottom left corner of~\figref{fig:Plate_40x24_kaSweep}. The pure-binary design obtained by the hard thresholding is included in the bottom right corner of~\figref{fig:Plate_40x24_kaSweep}. The consequences of the hard thresholding, \eg{}, performance drop, and resonance shift~\cite{diaz2010topology},~\cite{aage2010topology}, are moderated by the continuation scheme but not entirely suppressed. The binary design still exhibits low Q-factor~$\thr Q/{Q_\T{lb}^\T{TM}} = 1.06$, but the self-resonance is slightly shifted, see the frequency sweep in~\figref{fig:Plate_40x24_kaSweep}.

We repeatedly perform Q-factor minimization on the identical grid ($T=3840$) for $ka={0.1,0.4,0.8,1.2}$ to demonstrate the evolution of the optimized designs with electrical size. The designs are hard thresholded and shown in the first row of~\figref{fig:Other Designs}. It can be observed that a finer grid is required for low electric sizes~$ka = \{0.1,0.4\}$ to accommodate a structure with a self-resonant Q-factor, \eg{}, a meander line design~\cite{Capek_etal_2019_OptimalPlanarElectricDipoleAntennas}, to reach the lower bound. The optimizer realizes a self-resonant structure at~$ka=1.2$, but it perhaps requires a different feed position to get the lower bound. 

To demonstrate the effect of the density filter, we repeat the optimization for a fixed electric size of $ka=0.4$ while refining the discretization grid. Please refer to the thresholded designs in the second row of~\figref{fig:Other Designs}. The density filter with fixed radius~$\Rmin$ regularizes the solution space and can achieve structural convergence while the mesh is refined, see similarities in the designs. Nevertheless, the precise length-scale imposed by the density filter is lost while utilizing the continuation scheme of the projection filter~\cite{wang2011projection}. 

By visual inspection, the mesh with the finest granularity can accommodate a self-resonant meander line antenna; however, the density filter size must be further reduced, which may destabilize the convergence. 

\begin{figure}[!t]
\centering
\includegraphics[width=3.25in]{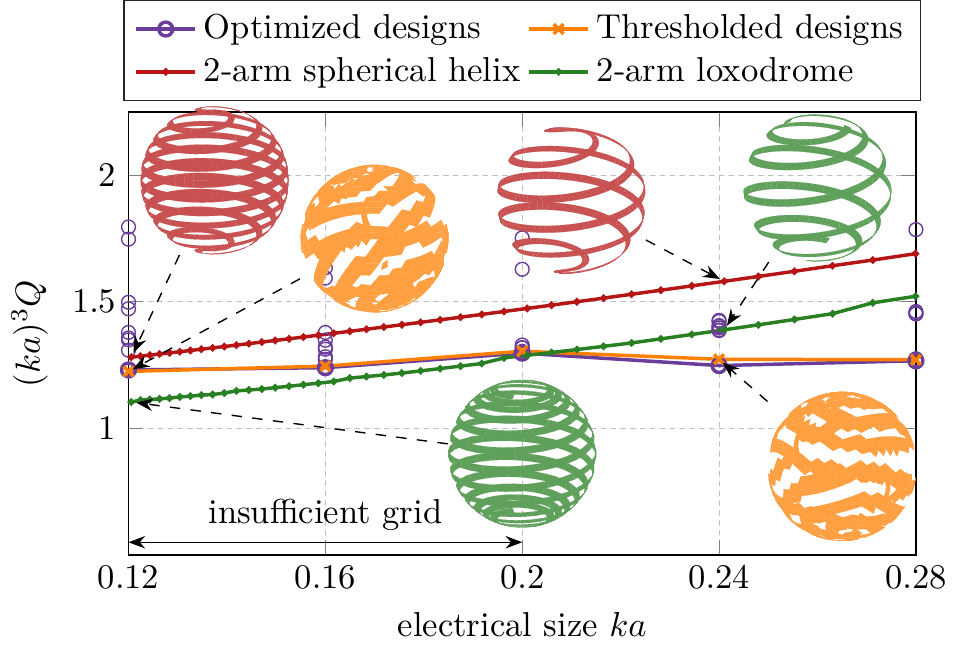}
\caption{Comparison of empirical and optimized designs, in terms of Q-factor values, are presented. The red line depicts self-resonant 2-arm spherical helix antennas, and the green line shows self-resonant 2-arm loxodrome antennas designed for a given electrical size~$ka$. Purple circles represent a single topology optimization run for the particular settings of~$\Rmin$ and~$\Sf$, \ie{}, ten runs are performed for each~$ka$. A purple line connects the best optimized designs, and the orange crosses depict the best thresholded structures. Notice that the discretization grid cannot accommodate a 2-arm geometry for low values of $ka$. For that reason, the topology optimization realizes a self-resonant 1-arm structure.}
\label{fig:Ka sweep designs}
\end{figure}

\subsection{Spherical Design Domain}

In the second example, the density topology optimization is utilized to solve the optimization problem~\eqref{eq:MMA task} for a spherical shell discretized into 2\,400~triangular elements (3\,600~basis function). The structure is excited by two identical delta gap sources placed symmetrically on the opposite sides of the equator to cover a shell uniformly with the material.
The performance of the optimized designs is assessed by the fundamental bound~\cite{CapekJelinek_OptimalCompositionOfModalCurrentsQ}, which for small electrical sizes reaches~$(ka)^3Q_\T{lb} \approx 1$. The optimization is performed with the maximal number of iterations~$I = 600$ for several electrical sizes in an interval~$ka\in[0.12,0.28]$, several radii of the density filter~$\Rmin=\{0.16,0.18,0.20,0.22\}a$, and two area fractions~$S_f=\{0.4,0.5\}$. As previously mentioned, it is often mandatory to restart the optimization with various initial parameters due to the non-convex nature of the problem.  

Optimized designs obtained by topology optimization (average computation time 50 minutes) for various parameters are shown in~\figref{fig:Ka sweep designs} by the purple circle marks. The purple line shows the performance of the best designs. Unsurprisingly, the optimizer finds a structure that combines two dominant spherical modes~\cite{CapekJelinek_OptimalCompositionOfModalCurrentsQ}. Thresholding the optimized designs results in a slight drop in performance, see orange line in~\figref{fig:Ka sweep designs}. For $ka<0.2$, the optimizer realizes self-resonant 1-arm designs and leaves one delta gap feeder isolated. A considerably finer mesh grid would be required to accommodate 2-arm geometry reaching the self-resonance. For $ka>0.2$, the discretization and the selected radius of the density filter form a solution space in which the proposed method finds self-resonant 2-arm spherical geometry with a low Q-factor value, \ie{}, for $ka=0.24$ $(ka)^3Q=1.25$. 

State-of-the-art spherical designs are helical antennas~\cite{Best_LowQelectricallySmallLinearAndEllipticalPolarizedSphericalDipoleAntennas}, see~\appref{app:parametrization} for their parametrization. Their Q-factors can further be reduced by folding, \ie{}, using a multi-arm spherical helix, which covers the spherical surface more uniformly\footnote{In~\cite{Best_LowQelectricallySmallLinearAndEllipticalPolarizedSphericalDipoleAntennas} the motivation for realizing the folded design was to increase the radiation resistance. Folding, nevertheless, also helps to approximate dominant spherical modes better.}. Compared with the topology-optimized designs, Q-factors of empirical 2-arm spherical helices for various electrical sizes~$ka$ are depicted in~\figref{fig:Ka sweep designs} by the red dots. 
The optimized designs outperform the empirical designs since the optimizer operates with more design variables. It is also worth noting that the optimized design does not resemble a spherical helix. The slope of a spherical helix is not constant along its curve, see~\appref{app:parametrization}. A curve conforming with the optimized design is loxodrome (Rhumb line)~\cite{alexander2004loxodromes} which always crosses the meridians at the same angle, see~\appref{app:parametrization} for the parametrization. {This agrees with the work~\cite{Kim_MinimumQESA} where it was analytically shown that an optimal current path for minimal Q-factor value follows loxodromes, but the authors improperly called such curves spherical helices. The authors of this manuscript recommend a differentiation between these curves since their Q-factor performance is not the same.} Loxodrome is only similar to a spherical helix near the equator but has an asymptotic point at the poles, \ie{}, the stereographic projection of a loxodrome to the plane of the equator is a logarithmic spiral. The Q-factor of the empirical 2-arm self-resonant loxodrome antennas for electrical sizes~$ka\in[0.12,0.28]$ are depicted in~\figref{fig:Ka sweep designs} by the green dots. It can be observed that loxodromes outperform the 2-arm spherical helix design in the entire $ka$ region. Curiously, they also outperform the optimized designs for $ka<0.2$. An insufficient mesh grid density used for the optimization causes this. For $ka>0.2$, the loxodrome has a higher Q-factor value than the optimized designs.

\begin{figure}[!t]
\centering
\includegraphics[width=3.25in]{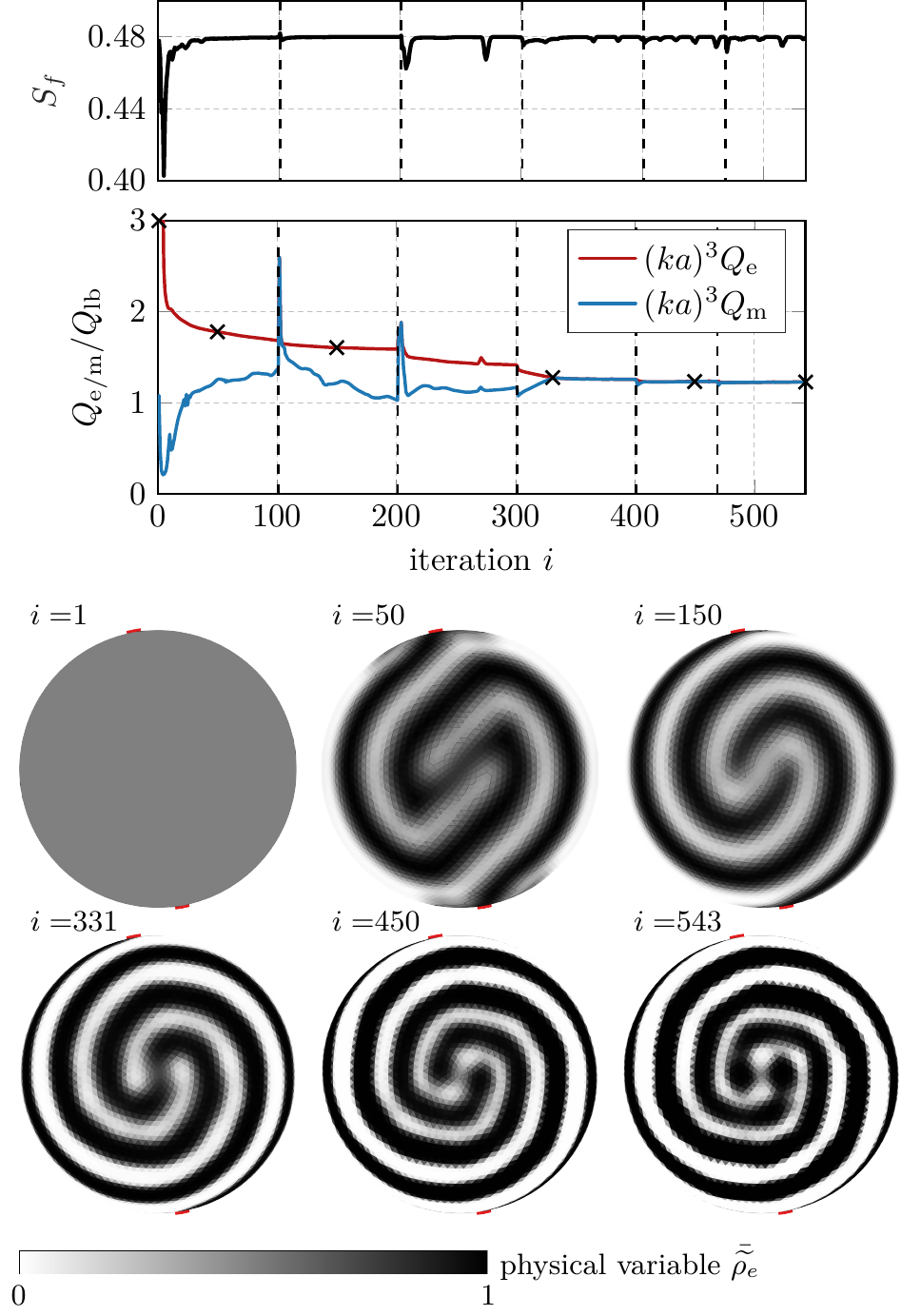}
\caption{(Top) Area fraction~$\Sf$, and (middle) convergence history of Q-factors $Q=\T{max}\{\Qe, \Qm\}$ normalized to the lower bound~$Q_\T{lb}$~\cite{CapekJelinek_OptimalCompositionOfModalCurrentsQ}~during the optimization. The dashed vertical lines represent iterations in which the sharpness~$\beta$ of projection filter~$\heavisideFilter$ was doubled. (Bottom) Snapshots of density field distribution from several iterations are included. The black line in the designs bounds the region where~$\desF\geq0.5$, thus showing a thresholded design boundary.}
\label{fig:Convergence_Helix}
\end{figure}

\begin{figure}[!t]
\centering
\includegraphics[width=3.25in]{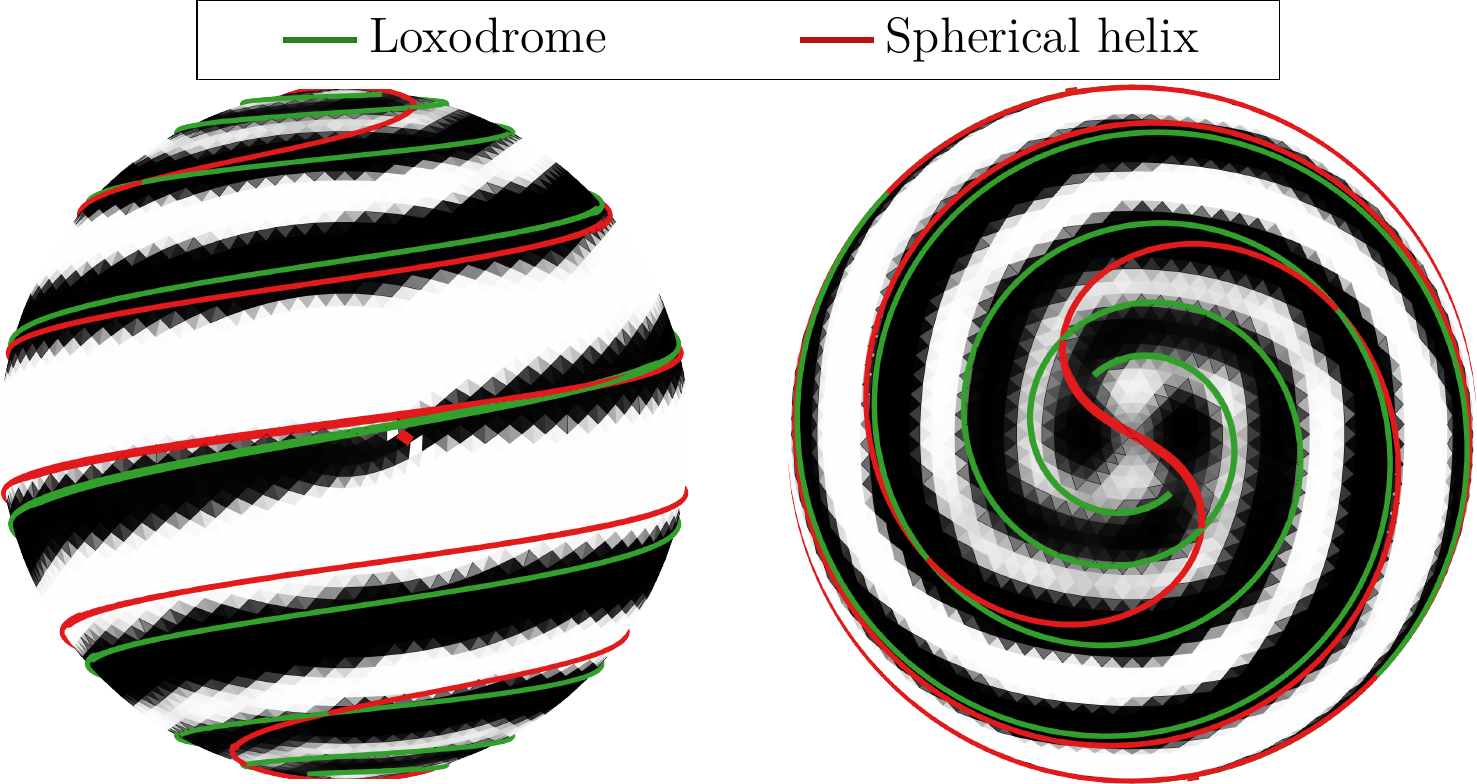}
\caption{The optimized design was obtained by density topology optimization. The structure's resemblance to the two spherical curves, \ie{}, spherical helix and loxodrome, is depicted. It is shown that the optimized topology follows the loxodrome curve more closely than the spherical helix. Both curves are very alike near the equator but distinct at the poles~{\cite{Kim_MinimumQESA}}.}
\label{fig:FineSphere}
\end{figure}

The optimization was repeatedly performed for electrical size $ka=0.24$ for a finer discretization of a spherical shell consisting of 11\,616 triangles (17\,424 basis functions) to obtain a higher level of structural details. The feeding configuration remains the same. The fixed radius of the density filter~$\Rmin=0.2a$ was used, and the optimization was subjected to the area fraction constraint\footnote{On the grid with $T=2400$ the optimizer distributed material using only 48\% of the surface, and thus, $\Sf$ is reduced to make the area fraction constraint active.}~$S_f=0.48$. The algorithm ran for 52 hours and successfully provided a design occupying 48\% of the design domain, see the first plot in~\figref{fig:Convergence_Helix}.  
The evolution of the optimized Q-factor is captured in the second plot in~\figref{fig:Convergence_Helix}, where~$\Qm$ and~$\Qe$ approach to realize a self-resonant topology.
The snapshots of the spherical designs during the optimization are included at the bottom of~\figref{fig:Convergence_Helix}. Within tens of iterations, the optimizer finds a 2-arm spherical topology whose arms are further extended by turning around the poles. The optimized structure has an optimized slope to combine dominant TM and TE spherical modes properly. The optimized design obtained by the density topology optimization exhibits~$(ka)^3Q=1.23$ and is shown in~\figref{fig:FineSphere} for two orthogonal views. The spherical helix and loxodrome curves are also plotted for comparison. 

{The analytical approach~\cite{Kim_MinimumQESA}} dictates an optimal current combination of dominant spherical modes, which forms a constant angle with the meridians. The loxodrome has the same properties and can approximate the optimal current distribution better than the spherical helix, see the difference of the curves mentioned above from the optimized topology in~\figref{fig:FineSphere}. The optimizer more accurately follows the loxodrome than the spherical helix. {The proposed method can deliver a close-to-optimal structure in terms of Q-factor value in an automated manner.}

\section{Conclusion}\label{sec:conclusion}
A complete treatment of density-based topology optimization in a method-of-moments framework was presented and verified on a particular problem of the Q-factor minimization utilizing conductor-based electrically small antennas. The proposed technique utilizes the projection operator in a continuation scheme to obtain a near-binary and standard density filter as an intermediate step to smooth an optimized design. The adjoint sensitivity analysis is exploited to attain gradient information to solve the optimization problem to local optimality. 

The efficacy of the presented approach was demonstrated by applying it to canonical examples of optimizing the Q-factor of rectangular and spherical design regions. The method excels in the scalability and efficient evaluation of gradients compared to other formulations of topology optimization. The realized performance is, in all cases, close to the fundamental bounds and self-resonance is met depending on the granularity of the discretization grid.
Furthermore, the proposed approach is capable of improving the performance of existing empirical designs in this case without modifying the nature of the topology. Alternatively, density topology optimization can inspire antenna engineers by providing novel structures used as starting points for future design methodology {or by proving the existence of well-performing design}. The framework can handle multi-port and multi-frequency optimization only at the cost of the linear increase of computational complexity.

The settings of the proposed method are tailored for Q-factor optimization, and may not work universally for a general electromagnetic metric. The optimization parameters significantly impact the optimization process and its success, and multiple runs with different settings are necessary. While the method has the advantage of fast convergence, the resulting designs are near-binary only, and a projection to a pure-binary result is required. Such a projection may result in a performance drop and resonance shift. Thus, the proposed technique can be conveniently combined into a two-step approach, in which the second step would tune the pure-binary design into self-resonance. 

There are many possibilities for future research that are closely connected to the proposed framework. For example, optimizing feeding position, amplitude, and phase can improve performance further and offer greater versatility. Other figures of merits should be explored as well. In the future, the framework could also be extended to include multiobjective optimization. Another potential direction is to define the topology optimization with a robust formulation~\cite{wang2011projection} so that it produces designs resilient to manufacturing errors. A final possibility is to extend the implementation of the method of moving asymptotes to cover the globally convergent version which is more effective for highly non-convex systems.

\appendices
\section{Electric Field Integral Equations within MoM formalism}\label{app:EFIE}
Assume a time-harmonic steady state at angular frequency~$\omega$ with convention~$\partial/\partial t \to \J\omega$, with $\J$ being the imaginary unit. Without loss of generality, surface formulation of EFIE~\cite{Gibson_MoMinElectromagnetics} for highly conducting bodies with space-dependent surface resistance $R_s(\V{r})$
is employed as
\begin{equation}
    -\normalvec\times \normalvec\times \left(\Ei(\V{r}) + \Es(\V{r})\right)  = R_s(\V{r}) \V{J}(\V{r}),
    \label{eq:EFIE}
\end{equation}
where $\normalvec$ is the surface normal unit vector, $\Ei(\V{r})$ is an incident field and $\Es(\V{r})$ is scattered field defined via
\begin{equation}
     \Es(\V{r})=-\J k Z_0 \int\limits_{\varOmega'} \overline{\M{G}}\left(\V{r},\V{r'}\right)\cdot\V{J}\left(\V{r'}\right)\D\V{r'}
     \label{eq:Escat}
\end{equation}
with $k$ and $Z_0$ being wavenumber and the free-space impedance, respectively. The symbol~$\overline{\M{G}}$ in \eqref{eq:Escat} denotes dyadic Green's function~\cite{Chew_WavevAndFieldsInInhomogeneousMedia} {with Cartesian components}

\begin{equation}
    {G}_{uv}\left(\V{r},\V{r'}\right) = \left[\delta_{uv} + \frac{1}{k^2}\dfrac{\partial^2}{\partial r_u \partial r_v}\right]\frac{\T{e}^{-\J k |\V{r}-\V{r'}|}}{4\pi|\V{r}-\V{r'}|},
\end{equation}

{where $\delta_{uv}$ is Kronecker delta and $r_u$ denotes the $u$-th component of radius vector~$\V{r}$.}

\begin{figure}[!t]
\centering
\includegraphics[width=3.25in]{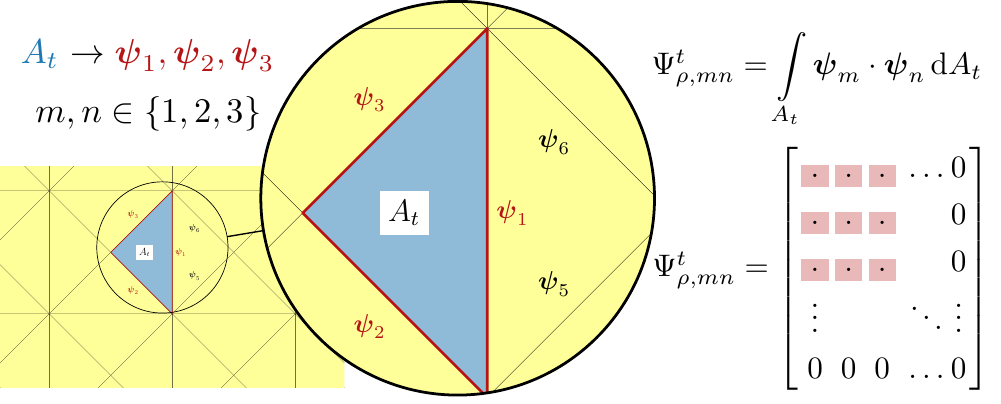}
\caption{The example of a particular material element matrix for triangle~$t$ is presented. This triangle is {shared} by three basis functions~$\{\basisFcn_1,\basisFcn_2,\basisFcn_3\}${, each of which also occupies one adjacent triangle.} {These basis functions have} unique labels {(given by the numbering of edges)} which determine certain positions with non-zero entries in element matrix~$\M{\Psi}_\rho^t$. The evaluation of the integral with RWG basis function can be found in~\cite{JelinekCapek_OptimalCurrentsOnArbitrarilyShapedSurfaces}.}
\label{fig:Lmat}
\end{figure}

The MoM process employs an expansion of unknown surface current density $\V{J}(\V{r})$ with a set of basis functions $\{\basisFcn_n(\V{r})\}$ into 
\begin{equation}
    \V{J}(\V{r}) \approx \displaystyle \sum_{n} I_n\basisFcn_n(\V{r})
\end{equation}
and Galerkin's testing procedure~\cite{Harrington_FieldComputationByMoM} to recast~\eqref{eq:EFIE} into its matrix form
\begin{equation}
    (\Zvac + \Zrho)\M{I} = \M{ZI} = \M{V},
    \label{appeq:ZI=V}
\end{equation}
where $\Zvac$ and $\Zrho$ are the vacuum and material parts of impedance matrix $\M{Z}$, respectively, and $\M{V}$ is the excitation vector. Explicitly, the operators in~\eqref{appeq:ZI=V} are defined element-wise as
\begin{align}
     Z_{0,mn} &= \J k Z_0 \int\limits_{\varOmega}\int\limits_{\varOmega'} \basisFcn_m(\V{r})\cdot \overline{\M{G}}\left(\V{r},\V{r'}\right) \cdot\basisFcn_n\left(\V{r'}\right)\D S'\D S,\\
     R_{\rho,mn} &= \displaystyle\int\limits_{\varOmega} R_s(\V{r}) \basisFcn_m(\V{r})\cdot\basisFcn_n(\V{r})\D S,\label{eq:Zrho}\\
     V_m &= \int\limits_{\varOmega}\basisFcn_m(\V{r})\cdot\V{E}^\text{i}\left(\V{r}\right)\D S. 
\end{align}
Without loss of generality, we utilize the triangular discretization and Rao-Wilton-Glisson basis functions~\cite{RaoWiltonGlisson_ElectromagneticScatteringBySurfacesOfArbitraryShape} {which are attached to each inner edge and its two adjacent triangles}. Within our design parametrization, the triangles have variable surface resistivity~$R_{\T{s}}(\des)$ according to the interpolation scheme~\eqref{eq:Interpolation}. Hence, we can rewrite~\eqref{eq:Zrho} as
\begin{equation}
\Zrho = \sum_{t=1}^T R_{\T{s}}(\des) \Zrhoe,
\end{equation}
where $\Zrhoe$ is a material element matrix that contains connectivity properties of a specific triangle~$t$. This matrix is written as
\begin{equation}
\Psi_{\rho,mn}^t = \int\limits_{A_t} \basisFcn_m\cdot\basisFcn_n\D A_t,
\label{eq:Zrhoe}
\end{equation}
where the indices~$m$ and~$n$ uniquely label the basis functions, forming a block matrix as shown in~\figref{fig:Lmat}. The integral in~\eqref{eq:Zrhoe} is evaluated analytically~\cite{JelinekCapek_OptimalCurrentsOnArbitrarilyShapedSurfaces}.

\section{Influence of the Interpolation scheme on Q-factor}\label{app:Choosing resistivities}

\begin{figure}[!t]
\centering
\includegraphics[width=3.25in]{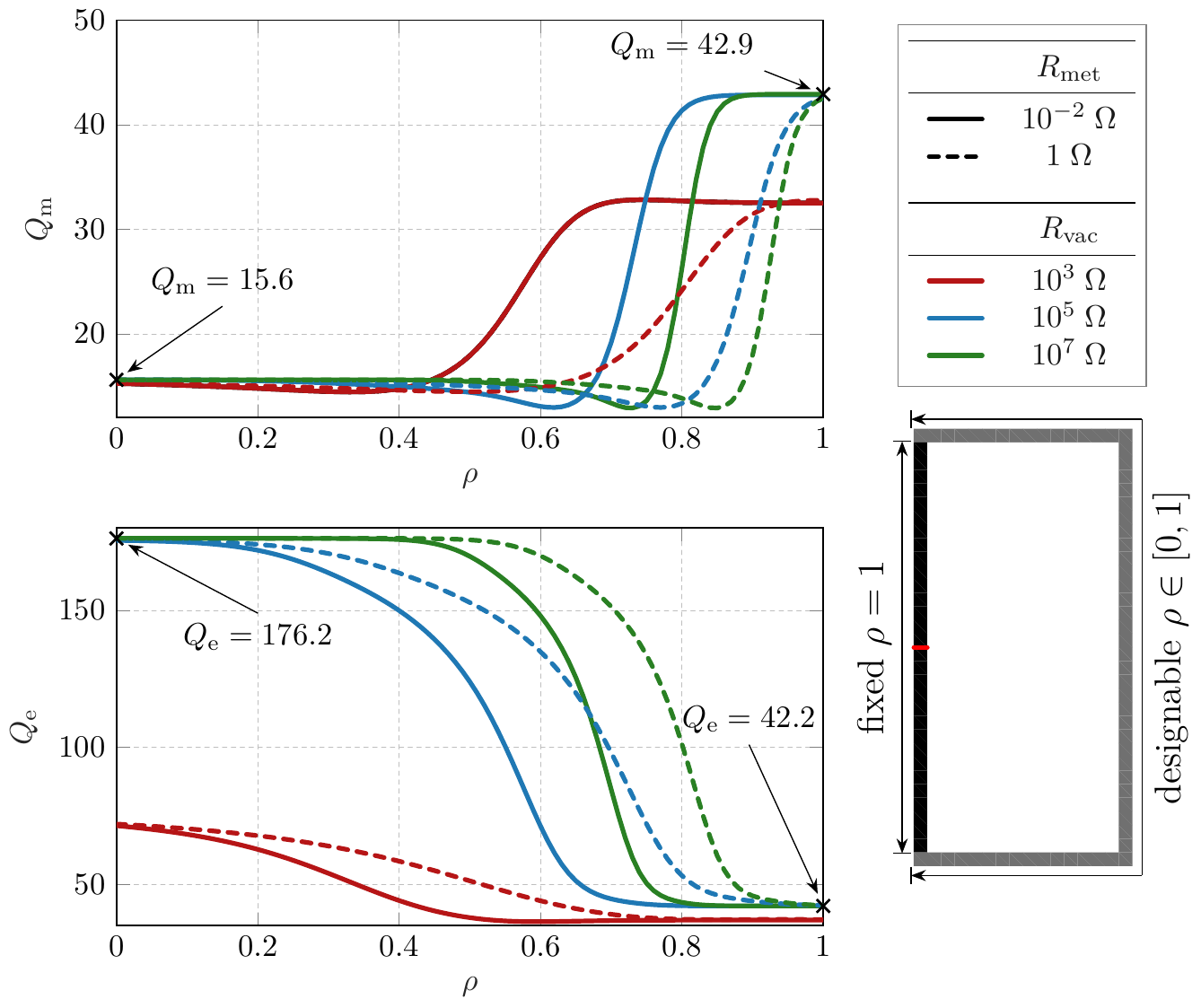}
\caption{The influence of interpolation scheme~\eqref{eq:Interpolation} on the electric and magnetic quality factor~$\Qem$ is examined at electrical size~$ka=0.59$. The loop antenna, self-resonant at the aforementioned~$ka$, is parametrized by variable~$\rho$ representing surface resistivity~$\Rs(\rho)$ through~\eqref{eq:Interpolation}. The structure is excited by a discrete delta gap placed in the middle of the left boundary. The antenna is divided into a fixed part with~$\rho=1$ and a designable part with~$\rho\in[0,1]$. The solid lines in the plot show the values of Q-factor~$\Qem$ for various combinations of the lower bounds~$\Rs(0)=\Zair$ as variable~$\rho$ is swept, with a fixed value of the upper bound~$\Rs(1) = \Zmet=10^{-2}\;\Omega$. Dashed lines depict the same sweep but with the upper bound value set to~$\Zmet=1\;\Omega$.} 
\label{fig:Behavior of Interpolation}
\end{figure}

The behavior of the interpolation scheme~\eqref{eq:Interpolation} on the Q-factor at electrical size~$ka=0.59$ is investigated using a model of a loop antenna. Physically accurate and admissible scenarios are discussed first. Assuming a loop antenna made of PEC close to resonance, \ie{}, $\Qe\approx\Qm$, the black crosses at~$\rho=1$ in~\figref{fig:Behavior of Interpolation} represent this scenario. Furthermore, by considering only the black portion of the structure in~\figref{fig:Behavior of Interpolation} made of PEC, with the remainder being a vacuum, an electrically short dipole is formed, and the black crosses indicate the corresponding Q-factor values at~$\rho=0$ in~\figref{fig:Behavior of Interpolation}.

The effect of the interpolation scheme~\eqref{eq:Interpolation} can be presented by a continuous change of surface resistivity~$\Rs(\rho)$ in the gray part of~\figref{fig:Behavior of Interpolation}. Q-factor~$\Qem$ is evaluated for various values of the lower bound~$\Rs(\rho=0)=\Zair$ and upper bound~$\Rs(\rho=1)=\Zmet$ on the attainable resistivities, as shown in~\figref{fig:Behavior of Interpolation}. It can be observed that the physically admissible choice of boundary resistivities is represented by the green and blue lines in~\figref{fig:Behavior of Interpolation}, which reliably reflect the physical behavior for~$\rho=1$ and~$\rho=0$. On the other hand, the combination represented by the red lines does not accurately reflect the physics at all. The blue and green lines differ mainly in the intermediate region of~$\rho$, where a small change in~$\rho$ should produce a measurable effect on the Q-factor evaluation. The blue dashed line, \ie{}, a combination of $[\Zair,\Zmet]=[10^5,1]\;\Omega$, delivers the most considerable freedom in the design variable for the optimizer while minimizing the Q-factor.

Selecting $[\Zair,\Zmet]=[10^5,1]\;\Omega$ as boundary resistivity values results in an inaccurate modeling of the physical behavior. However, based on the authors' observations while analyzing an optimized design with proper resistivity values, such as $[\Zair,\Zmet]\approx[\infty,0]\;\Omega$, there is only a negligible drop in Q-factor performance. Therefore, it is advantageous to ease the optimization process at the expense of the small accuracy loss in Q-factor evaluation. It is important to stress that this claim is not universal and depends on the interpolation scheme or optimized metric. Hence, analogous investigations are always recommended.

\section{Q-Factor and its Adjoint Sensitivi
ty Analysis}\label{app:Quality factor}
The {radiation} Q-factor of an antenna can be defined according to~\cite{CollinRotchild_EvaluationOfAntennaQ} as 
\begin{equation}
    Q = \frac{2 \omega \max\{\We,\Wm\}}{\Prad} \equiv \max\{\Qe,\Qm\},
    \label{eq:Q}
\end{equation}
where $\omega$ is angular frequency, $\We$ and $\Wm$ are cycle mean stored electric and magnetic energy, respectively, and $\Prad$ is the cycle mean radiated power
\begin{equation}
    \Prad = \frac{1}{2}\M{I}^\herm \M{R}_0 \M{I},
\end{equation}
where $\M{R}_0$ is a real part of the vacuum impedance matrix and superscript $^\herm$ is the complex conjugate transpose. 
Stored energies in~\eqref{eq:Q} are approximated~\cite{Vandenbosch_ReactiveEnergiesImpedanceAndQFactorOfRadiatingStructures},~\cite{Gustafsson_OptimalAntennaCurrentsForQsuperdirectivityAndRP} as
\begin{align}
    \We &\approx \frac{1}{8}\M{I}^\herm\left(\frac{\partial\M{X}_0}{\partial \omega} - \frac{\M{X}_0}{\omega}\right)\M{I} = \frac{1}{4\omega} \M{I}^\herm \XEmat \M{I},\\
    \Wm &\approx \frac{1}{8}\M{I}^\herm\left(\frac{\partial\M{X}_0}{\partial \omega} + \frac{\M{X}_0}{\omega}\right)\M{I} = \frac{1}{4\omega} \M{I}^\herm \XMmat \M{I},
\end{align}
$\M{X}_0$ is an imaginary part of the vacuum impedance matrix, and $\M{X}_0 = \XMmat - \XEmat$. These matrix operators are obtained from method of moments codes~\cite{atom}. The definition~\eqref{eq:Q} is interpreted as the upper bound on the Q-factor of an antenna system that is tuned to the resonance by the addition of a~lossless reactive element. 

The sensitivities of $\Qe$ and $\Qm$ with respect to varible~$\des$ can be computed by adjoint sensitivity analysis~\cite{tortorelli1994design} as
\begin{equation}
        \frac{\D \Qem}{\D \des} =\frac{\partial \Qem}{\partial \des} + 2 \RE \left\{\V{\lambda}_\T{e/m}^\trans\frac{\partial \Zrho}{\partial \des}\M{I} \right\},
        \label{eq:sensitivities appendix}
\end{equation}
where $\adj_\T{e/m}$ is the solution to corresponding adjoint equations
\begin{equation}
    \M{Z}^\trans\adj_\T{e/m} = -\left(\frac{\partial \Qem}{\partial \M{I}}\right)^\trans,
\end{equation}
where the right-hand side reads
\begin{equation}
\frac{\partial \Qem}{\partial \M{I}} = \frac{1}{2\Prad}\M{I}^\herm\left(\M{X}_\T{e/m} - \Qem\M{R}_0 \right).
\end{equation}
The Q-factor is not a function of~$\des$, and the first term in~\eqref{eq:sensitivities appendix} vanishes.

\section{{Resilience towards initial density distribution}}\label{app:initial seed}
\begin{figure}[!t]
\centering
\includegraphics[width=3.25in]{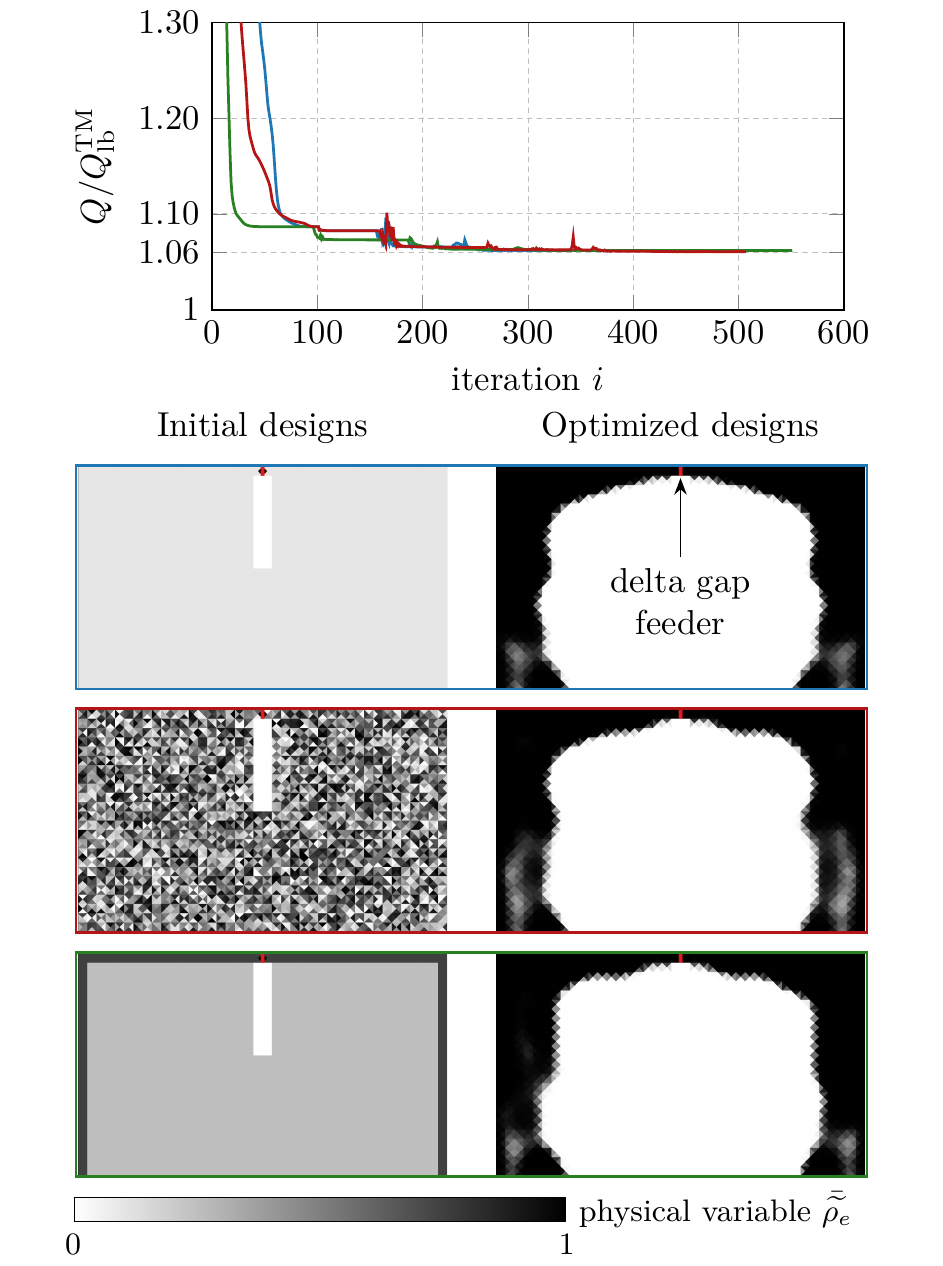}
\caption{{(Top) Convergence history of Q-factors $Q=\T{max}\{\Qe, \Qm\}$ normalized to the fundamental bound~$Q_\T{lb}^{\T{TM}}$~\cite{GustafssonTayliEhrenborgEtAl_AntennaCurrentOptimizationUsingMatlabAndCVX} for three different initial density distribution. Occasional jumps result from the projection filter update. (Bottom) Snapshots of an initial and optimized density field distribution are included below in the convergence graph. The outline of each design is color-matched with the curve in the convergence plot.}} 
\label{fig:convergence initial seed}
\end{figure}

{The Q-factor minimization by the density-based topology optimization on a rectangular design region from~\secref{sec:rectangle} is considered here with the same settings but different initial density distribution, namely, uniform distribution with~$\des=0.1$, random distribution, and distribution with already well-separated charge distribution (low Q-factor value), see the bottom of~\figref{fig:convergence initial seed}. The convergence curves are reported at the top of~\figref{fig:convergence initial seed}. The method delivered optimized designs, see bottom of~\figref{fig:convergence initial seed}, with similar values of Q-factor~$Q/Q_\T{lb}^\T{TM}=1.06$, \ie{}, the nature of Q-factor and usage of spatial filtering make the proposed method less dependent on the initial seed. The main influence of the initial seed lies in the convergence sharpness at the early iterations. This study might nonetheless be not applicable to EM metrics different from the Q-factor.   
}

\section{Parametrization of Spherical Curves}\label{app:parametrization}
A spherical helix~\cite{scofield1995curves} is a curve maintaining a constant angle with the equatorial plane and is described by the number of turns~$M$ which connect the poles. A parametrization of a spherical helix is given by 
\begin{equation}
\label{eq:app:spherHelix}
    \V{r}(t) = R \Big[f(t)\cos(2\pi Mt), f(t)\sin(2\pi Mt),2t-1 \Big],
\end{equation}
where~$f(t)=\sqrt{1-4(t-0.5)^2}$, $R$ is the radius of the sphere and where different components of~\eqref{eq:app:spherHelix} refer to Cartesian coordinates.

Loxodrome (Rhumb line)~\cite{alexander2004loxodromes} is a curve that crosses all meridians at the same angle. Its parametrization in {Cartesian coordinates} reads 
\begin{equation}
    \V{r}(t) = R\left[\dfrac{\cos(t)}{\cosh(\gamma t)}, \dfrac{\sin(t)}{\cosh(\gamma t)}, \tanh(\gamma t)\right],
\end{equation}
{where~$\gamma$ is the slope of the curve. Alternative parametrization in spherical coordinates is given in~\cite{Kim_MinimumQESA}.}
Denoting~$\V{\tau} = \T{d} \V{r} / \T{d} t $ a tangent vector field along the curve, it can be shown that
\begin{equation}
    \left| \V{\tau} \right| = \dfrac{\sqrt{1+\gamma^2}}{\cosh \left( \gamma t\right)},
\end{equation}
and that spherical vector components of the unit vector along the curve read 
\begin{equation}
    \dfrac{\tau_r}{\left| \V{\tau} \right|} = 0, \, \dfrac{\tau_\theta}{\left| \V{\tau} \right|} = -\dfrac{\gamma}{\sqrt{1+\gamma^2}}, \,  \dfrac{\tau_\varphi}{\left| \V{\tau} \right|} = \dfrac{1}{\sqrt{1+\gamma^2}}.
\end{equation}
This is the vector orientation of an addition of~$\T{TM}_{10}$ and ~$\T{TE}_{10}$ spherical vector waves, \ie{}., the correct orientation to minimize Q-factor on a spherical shell~{\cite{Kim_MinimumQESA}},~\cite{CapekJelinek_OptimalCompositionOfModalCurrentsQ}.

\ifCLASSOPTIONcaptionsoff
  \newpage
\fi


\begin{IEEEbiography}
[{\includegraphics[width=1in,height=1.25in,clip,keepaspectratio]{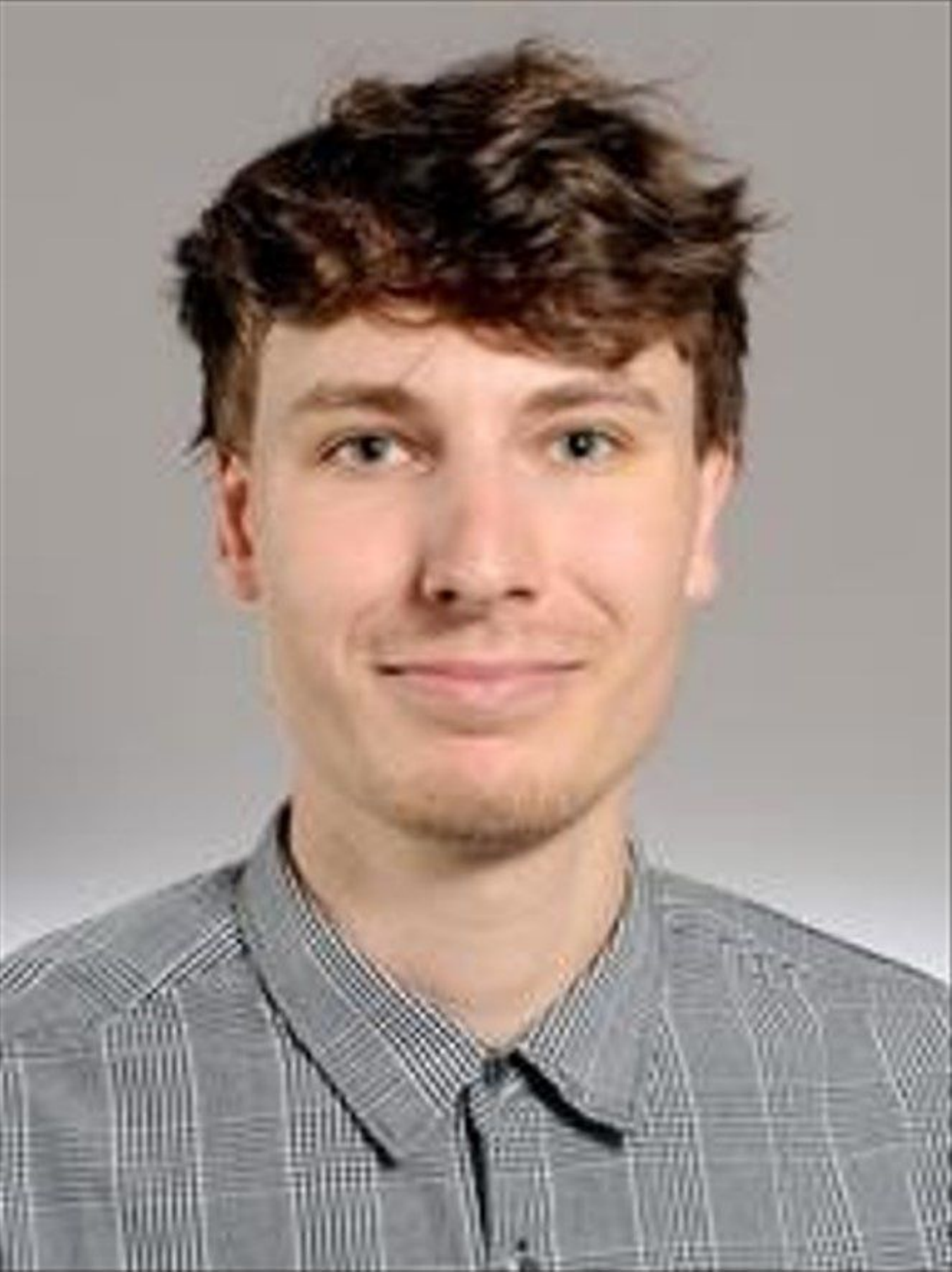}}]{Jonas Tucek}
received the M.Sc. degree in electrical engineering from Czech Technical University, Prague, Czech Republic, in 2021, where he is currently pursuing the Ph.D. degree in the area of topology optimization. 
\end{IEEEbiography}

\begin{IEEEbiography}[{\includegraphics[width=1in,height=1.25in,clip,keepaspectratio]{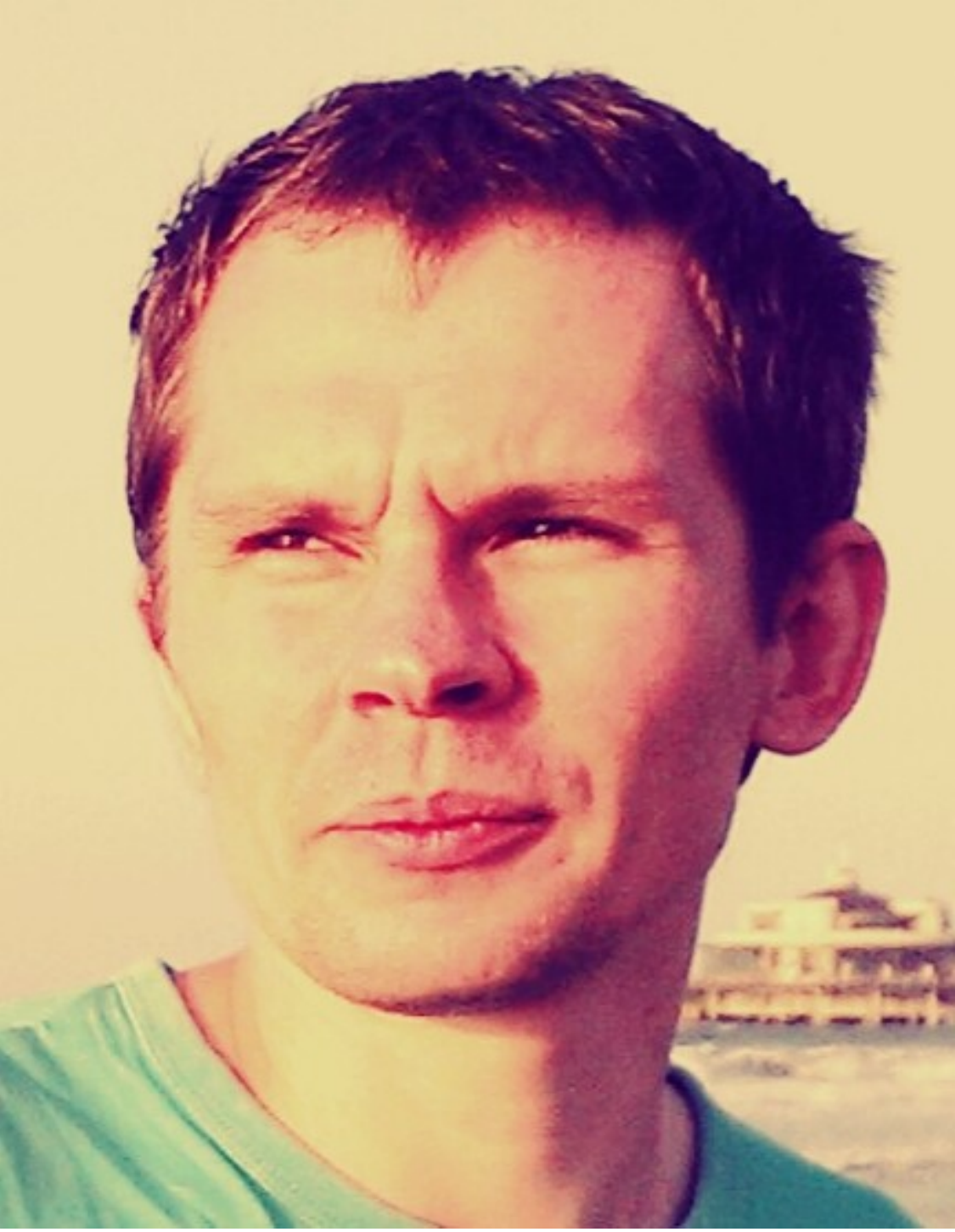}}]{Miloslav Capek}
(M'14, SM'17) received the M.Sc. degree in Electrical Engineering 2009, the Ph.D. degree in 2014, and was appointed Associate Professor in 2017, all from the Czech Technical University in Prague, Czech Republic.
	
He leads the development of the AToM (Antenna Toolbox for Matlab) package. His research interests are in the area of electromagnetic theory, electrically small antennas, antenna design, numerical techniques, and optimization. He authored or co-authored over 130~journal and conference papers.

Dr. Capek is the Associate Editor of IET Microwaves, Antennas \& Propagation. He was a regional delegate of EurAAP between 2015 and 2020 and an associate editor of Radioengineering between 2015 and 2018.
\end{IEEEbiography}

\begin{IEEEbiography}[{\includegraphics[width=1in,height=1.25in,clip,keepaspectratio]{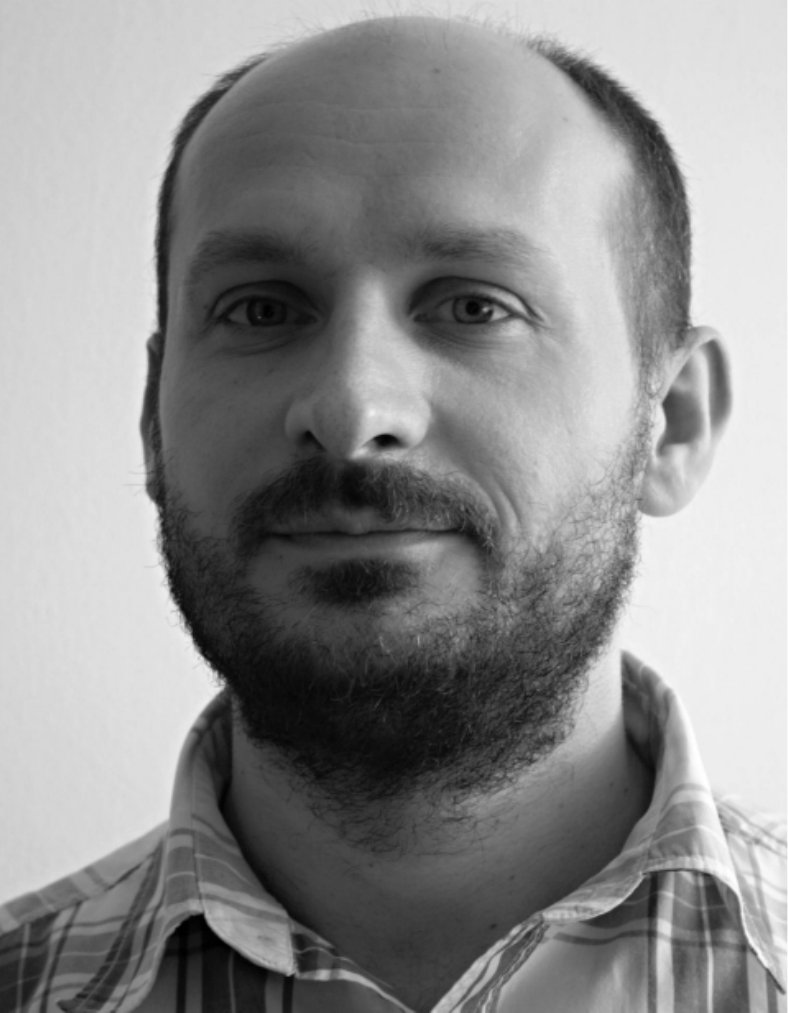}}]{Lukas Jelinek}
received his Ph.D. degree from the Czech Technical University in Prague, Czech Republic, in 2006. In 2015 he was appointed Associate Professor at the Department of Electromagnetic Field at the same university.

His research interests include wave propagation in complex media, electromagnetic field theory, metamaterials, numerical techniques, and optimization.
\end{IEEEbiography}


\begin{IEEEbiography}
[{\includegraphics[width=1in,height=1.25in,clip,keepaspectratio]{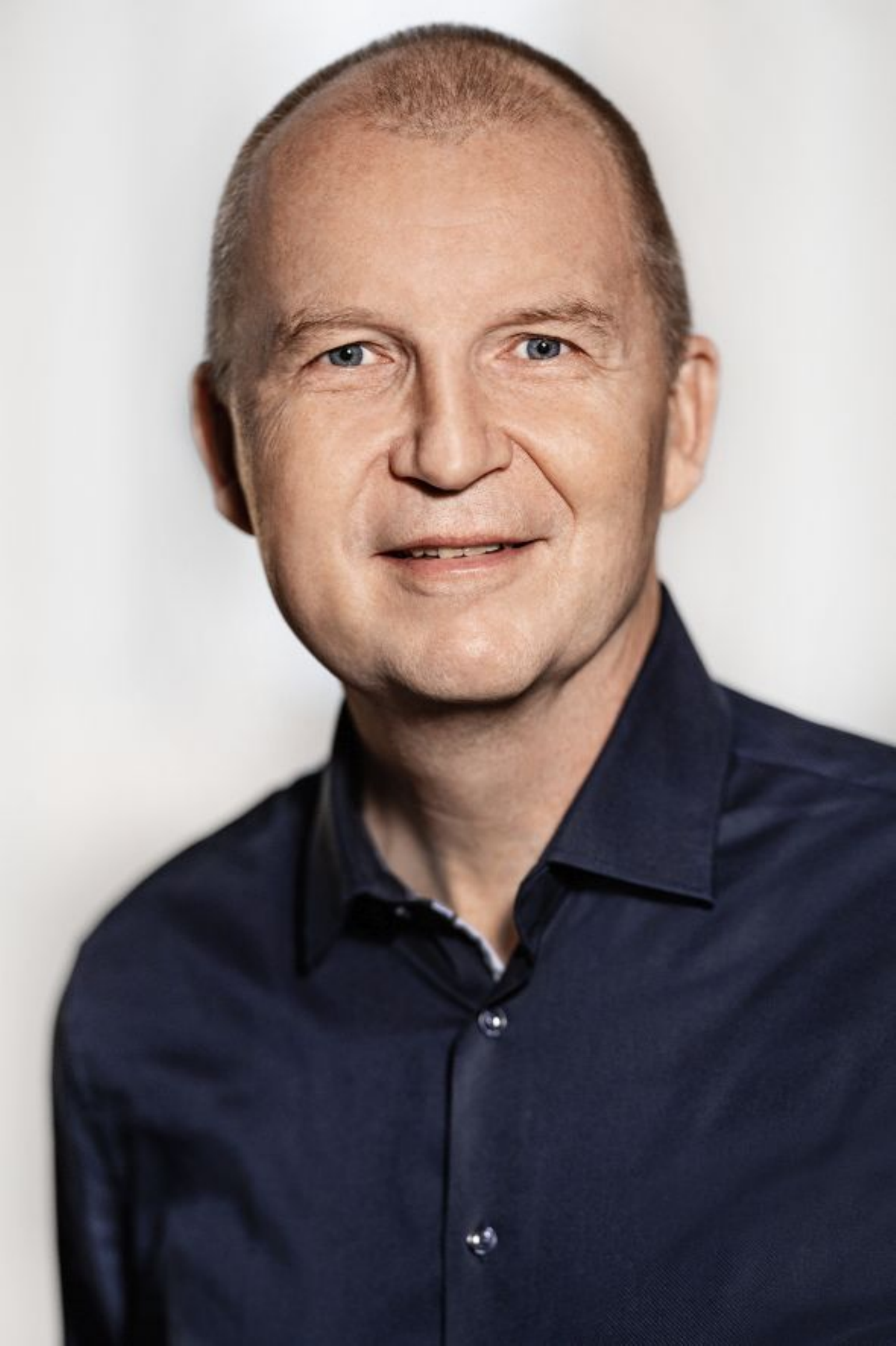}}]{Ole Sigmund}
 received the Ph.D. degree and Habilitation in 1994 and 2001, respectively, and has had research positions with the University of Essen and Princeton University. He is a professor in the Department of Mechanical Engineering, Technical University of Denmark (DTU). He is a member of the Danish Academy of Technical Sciences and the Royal Academy of Science and Letters, Denmark and is the former elected President (2011-15, now EC member) of ISSMO (International Society of Structural and Multidisciplinary Optimization). 

His research interests include theoretical extensions and applications of topology optimization methods to mechanics, and multiphysics problems.
\end{IEEEbiography}
\end{document}